\documentclass[11pt]{article}
\usepackage[margin=1in]{geometry}

\usepackage{authblk} 
\usepackage{graphicx}
\usepackage{amssymb}

\usepackage{amsmath}
\usepackage{enumerate}
\usepackage{amsfonts}
\usepackage[mathcal]{euscript}
\usepackage{amsthm}
\usepackage[all,2cell,import]{xy} \UseAllTwocells
\usepackage{xcolor}
\usepackage{tkz-euclide}
\usepackage{tikz}
	\usetikzlibrary{3d} 
	\usetikzlibrary{decorations.pathreplacing,angles,quotes}

\usepackage{url}
\makeatletter
\def\url@leostyle{%
	\@ifundefined{selectfont}{\def\UrlFont{\sf}}{\def\UrlFont{\small\ttfamily}}}
\makeatother
\usepackage[refpage]{nomencl}

\makenomenclature
\usepackage{makeidx}
\makeindex
\numberwithin{equation}{section}

\usepackage[colorlinks=true, pdfstartview=FitV, linkcolor=black,
			   citecolor=black, urlcolor=black]{hyperref}
\usepackage{bookmark}
\usepackage[titletoc]{appendix} 
\hypersetup{hypertexnames=false} 

          \newcommand{\nc}{\newcommand}
          \nc{\DMO}{\DeclareMathOperator}	
          \nc{\commentout}[1]{}

          \nc{\newnotation}{\nomenclature}
          
          \nc{\Diff}{\mathcal{D}\!\operatorname{iff}}
          \DMO{\emb}{Emb}
          \nc{\Ind}{\mathrm{Ind}}
          \nc{\Loc}{\mathsf{Loc}}
          \nc{\Cob}{\mathsf{Cob}}
          \nc{\mul}{\mathsf{Mul}}
          \nc{\fat}{\mathsf{fat}}
          \nc{\cob}{\mathsf{Cob}}
          \nc{\coh}{\mathsf{Coh}}
          \nc{\idem}{\mathsf{Idem}}
          \nc{\sets}{\mathsf{Sets}}
          \nc{\near}{\mathsf{near}}
          \nc{\sing}{\mathsf{Sing}}
          \nc{\Sing}{\mathsf{Sing}}
          \nc{\perf}{\mathsf{Perf}}
          \nc{\block}{\mathsf{block}}
          \nc{\ssets}{\mathsf{sSets}}
          \nc{\cmpct}{\mathsf{cmpct}}
          \nc{\compact}{\mathsf{cmpct}}
          \nc{\pwrap}{\mathsf{PWrap}}
          \nc{\coder}{\mathsf{Coder}}
          \nc{\bimod}{\mathsf{Bimod}}
          \nc{\grmod}{\mathsf{GrMod}}
          \nc{\Morita}{\mathsf{Morita}}
          \nc{\morita}{\mathsf{Morita}}
          \nc{\spaces}{\mathsf{Spaces}}
          \nc{\posets}{\mathsf{Poset}}
                    \nc{\Fun}{\mathsf{Fun}}
          \nc{\fun}{\mathsf{Fun}}
          \nc{\vect}{\mathsf{Vect}}
          \nc{\chain}{\mathsf{Chain}}
          \nc{\chainn}{Chain}
          \nc{\wrfuk}{\mathsf{WrFukaya}}
          \nc{\wrfukcompact}{\mathsf{WrFukaya}_{\mathsf{cmpct}}}
          \nc{\pwrfuk}{\mathsf{PWrFukaya}}
          \nc{\inffuk}{\mathsf{InfFuk}}
          \nc{\pwrfukml}{\mathsf{PWrFukaya}_{M,\Lambda}}
          \nc{\inffukml}{\mathsf{InfFuk}_{M,\Lambda}}
          \nc{\nattrans}{\mathsf{NatTrans}}
          \nc{\corres}{\mathsf{Corres}}

          \nc{\cat}{\mathsf{Cat}}
          \nc{\Cat}{\mathsf{Cat}}
          \nc{\ainfty}{\mathsf{A}_\infty}
          \nc{\inftycat}{\mathcal{C}\!\operatorname{at}_\infty}
          \nc{\inftyCat}{\mathcal{C}\!\operatorname{at}_\infty}
          \nc{\inftyGpd}{\mathcal{G}\!\operatorname{pd}_\infty}
          \nc{\Ainftycat}{\mathcal{C}\!\operatorname{at}_{A_\infty}}
          \nc{\dgcat}{\mathcal{C}\!\operatorname{at}_{dg}}
          \nc{\ainftycat}{\mathcal{C}\!\operatorname{at}_{A_\infty}}
          \nc{\stablecat}{\mathcal{C}\!\operatorname{at}_\infty^{\Ex}}

          \DMO{\op}{op}
          \DMO{\sd}{sd} 
          \DMO{\im}{im}
          \DMO{\ev}{ev}
          \DMO{\st}{st}
          \DMO{\stable}{Ex}
          \DMO{\map}{Map}
          \nc{\sdcoll}{\sd^{\coll}}
          \nc{\Excoll}{\Ex^{\coll}}
          \DMO{\inj}{inj}
          \DMO{\fib}{fib}
          \DMO{\conf}{Conf}
          \DMO{\chains}{Chains}
          \DMO{\cochains}{Cochains}
          \DMO{\cone}{Cone}
          \DMO{\Map}{Map}
          \DMO{\ran}{Ran}
          \DMO{\rot}{Rot}
          \DMO{\leg}{Leg}
          \DMO{\imm}{Imm}
          \DMO{\adj}{adj}
          \DMO{\symp}{Symp}
          \DMO{\tree}{Tree}
          \DMO{\cube}{Cube}
          \DMO{\weak}{weak}
          \DMO{\strong}{strong}
          \DMO{\Hoch}{Hoch}
          \DMO{\front}{front}
          \DMO{\flow}{Flow}
          \DMO{\floer}{Floer}
          \DMO{\Maps}{Maps}
          \DMO{\exact}{exact}
          \DMO{\excess}{Excess}
          \DMO{\Decomp}{Decomp}
          \DMO{\decomp}{Decomp}
          \DMO{\collar}{collar}
          \DMO{\yoneda}{Yoneda}
          \DMO{\hamspace}{Ham}
          \DMO{\sympspace}{Symp}
          \DMO{\holomaps}{Holomaps}
          \DMO{\comp}{Comp}
          \DMO{\crit}{Crit}
          \DMO{\test}{{test}}
          \DMO{\sign}{sign}
          \DMO{\topp}{top}
          \DMO{\indx}{Index}
          \DMO{\Break}{Break} 
          \DMO{\zero}{zero} 
          \DMO{\ob}{Ob}
          \DMO{\gr}{Gr} 
          \DMO{\Gr}{Gr} 
          \DMO{\cl}{Cl} 
          \DMO{\grlag}{GrLag}
          \DMO{\Pin}{Pin}
          \DMO{\Graph}{Graph}
          \DMO{\pin}{Pin}
          \DMO{\gap}{Gap}
          \DMO{\Ex}{Ex}
          \DMO{\id}{id}
          \DMO{\End}{End}
          \DMO{\sym}{Sym}
          \DMO{\aut}{Aut}
          \DMO{\Aut}{Aut}
          \DMO{\fine}{fine}
          \DMO{\haut}{hAut}
          \DMO{\hAut}{hAut}
          \DMO{\DK}{DK} 
          \DMO{\poly}{poly} 
          \DMO{\diff}{Diff}
          \DMO{\coll}{coll}
          \DMO{\dist}{dist} 
          \DMO{\coker}{coker} 
          \nc{\kernel}{\ker} 
          \DMO{\sspan}{span}
          \DMO{\hocolim}{hocolim}	
          \DMO{\holim}{holim}
          \DMO{\sk}{sk}

          \DMO{\ho}{ho}
          \DMO{\fin}{fin}
          \DMO{\tor}{Tor}
          \DMO{\ext}{Ext}
          \DMO{\ret}{Ret}
          \DMO{\ham}{Ham}
          \DMO{\con}{con}
          \DMO{\leaf}{leaf}
          \DMO{\supp}{supp}
          \DMO{\edge}{edge}
          \DMO{\colim}{colim}
          \DMO{\edges}{edges}
          \DMO{\Image}{image}
          \DMO{\roots}{roots}
          \DMO{\height}{height}
          \DMO{\finmod}{FinMod}
          \DMO{\leaves}{leaves}
          \DMO{\planar}{planar}
          \DMO{\vertices}{vertices}

          \nc{\lagg}{\lag^{\cG}}
          \nc{\iso}{\mathsf{Iso}}
          \nc{\Set}{\mathsf{Set}}
          \nc{\Ass}{\mathsf{ \bf Ass}}
          \nc{\Mod}{\mathsf{Mod}}
          \nc{\modules}{\mathsf{Mod}}
          \nc{\sset}{\mathsf{sSet}}
          \nc{\liou}{\mathsf{Liou}}
          \nc{\poset}{\mathsf{Poset}}
          \nc{\trno}{T^*\RR^n_{\geq 0}}
          \nc{\spectra}{\mathsf{Spectra}}
          \nc{\tensorfin}{\tensor^{\fin}}
          \nc{\lagptg}{\lag_{pt,pt}^{\cG}}
          \nc{\Fin}{\mathcal{F}\mathsf{in}}
          \nc{\lagnl}{\lag_{N,\Lambda}}
          \nc{\lagmlg}{\lag_{M,\Lambda}^{\cG}}
          \nc{\lagsplit}{\lag^{\mathsf{split}}}
          \nc{\lagktimes}{(\lag^{\dd k})^\times}
          \nc{\lagplanar}{\lag^{\times,\planar}}
          \nc{\Cont}{\text{\rm Cont}}
          \nc{\Ham}{\text{\rm Ham}}
          \nc{\Dev}{\text{\rm Dev}}
          \nc{\Lin}{\text{\rm Lin}}
          \nc{\Int}{\text{\rm Int}}
          \nc{\Hom}{\text{\rm Hom}}
          \nc{\Chord}{\text{\rm Chord}}
          \nc{\nbhd}{\mathcal{N}\text{\rm{bhd}}}

          \nc{\smsh}{\wedge}
          \nc{\un}{\underline}
          \nc{\xto}{\xrightarrow}
          \nc{\xra}{\xto}
          \nc{\tensor}{\otimes}
          \nc{\del}{\partial}
          \nc{\dd}{\diamond}
          \nc{\tri}{\triangle}
          \nc{\bb}{\Box}
          \nc{\into}{\hookrightarrow}
          \nc{\onto}{\twoheadrightarrow}
          \nc{\contains}{\supset}
          \nc{\transverse}{\pitchfork}
          \nc{\uncirc}{\underline{\circ}}

          \nc{\hiro}{\textcolor{blue}}

          \nc{\eqn}{\begin{equation}}
          \nc{\eqnn}{\begin{equation}\nonumber}
          \nc{\eqnd}{\end{equation}}
          \nc{\enum}{\begin{enumerate}}
          \nc{\enumd}{\end{enumerate}}
          \nc{\beastar}{\begin{eqnarray*}}
          \nc{\eeastar}{\end{eqnarray*}}

          \def\cC{\mathcal C}\def\cD{\mathcal D}
          \def\cG{\mathcal G}

          \def\AA{\mathbb A}\def\BB{\mathbb B}\def\CC{\mathbb C}
          \def\HH{\mathbb H}

          \def\RR{\mathbb R}

          \nc{\Euc}{\mathsf{Euc}}
          \nc{\DTop}{\mathsf{DTop}}
          \nc{\simp}{\mathsf{Simp}}
          \nc{\Ainftycatt}{A_\infty Cat}
          \nc{\dgcatt}{dg Cat}
          \nc{\StableCat}{StableCat}
          \nc{\subdivision}{\mathsf{subdiv}}
          \nc{\Kan}{\mathcal{K}\mathsf{an}}

          \nc{\deR}{\operatorname{deR}}

          \theoremstyle{definition}
          \newtheorem{theorem}{Theorem}[section]
          
          \newtheorem{prop}[theorem]{Proposition}
          \newtheorem{lemma}[theorem]{Lemma}

          \newtheorem{corollary}[theorem]{Corollary}
          \newtheorem{construction}[theorem]{Construction}

          \newtheorem{defn}[theorem]{Definition}
          \newtheorem{notation}[theorem]{Notation}
          \newtheorem{example}[theorem]{Example}
          \newtheorem{choice}[theorem]{Choice}
          \newtheorem{setup}[theorem]{Set-up}

          \newtheorem{remark}[theorem]{Remark}

\DMO{\fr}{Fr}
\DMO{\img}{Image}
\DMO{\KKan}{Kan}
\nc{\kan}{\mathsf{Kan}}
\nc{\quillen}{\mathsf{Quillen}}
\DMO{\Fr}{Fr}
\DMO{\finite}{finite}

\nc{\MMfld}{\mathsf{Mfld}}
\nc{\TTop}{\mathsf{Top}}

\nc{\mfld}{\mathcal{M}\!\operatorname{fld}}
\nc{\Mfld}{\mathcal{M}\!\operatorname{fld}}
\nc{\Top}{\mathcal{S}\!\operatorname{pc}}

\nc{\stab}{+1}
\nc{\stabn}{+n}

\nc{\embtop}{\mathcal{E}\kern-0.22em \operatorname{mb}}

\begin{document}

\title{Spaces over $BO$ are thickened manifolds}
\author{Hiro Lee Tanaka}
\maketitle

\begin{abstract}
Consider the topologically enriched category of compact smooth manifolds (possibly with corners), with morphisms given by codimension zero smooth embeddings. Now formally identify any object $X$ with its thickening $X \times [-1,1]$. 
We prove that the resulting $\infty$-category of thickened smooth manifolds is equivalent to the $\infty$-category of finite spaces over $BO$. (This is one formalization of the philosophy that embedding questions become homotopy-theoretic upon passage to higher dimensions.)  The central tool is a geometric construction of pushouts in this $\infty$-category, carried out with an eye toward proving analogous results in exact symplectic geometry. Notably, the proof never invokes smooth approximation nor any $h$-principle.
\end{abstract}

\tableofcontents

\clearpage

\section{Introduction}

Throughout this work, we adopt the convention that a smooth manifold may have boundary and corners, and that all connected components of a manifold have equal dimension.

Consider the $\infty$-category $\mfld^{\dd}$ whose objects are smooth, compact manifolds $X$, up to the relation 
	\eqn\label{eqn. thickening X}
	X \sim X \times [-1,1],
	\eqnd 
and where morphisms are codimension zero smooth embeddings---i.e., smooth embeddings whose derivatives point-wise have trivial cokernel. We do not require that the maps respect boundary or corner strata in any way. See Notation~\ref{notation. thickened manifold category} and~\ref{notation. mfld thickened} for a model of this $\infty$-category.

\begin{remark}
\label{remark. points are disks}
Let $X = [-1,1]^n$ be the cube and $Y = D^n $ the closed unit disk. Appropriate scalings define a smooth codimension zero embedding of $X$ inside $Y$, and of $Y$ inside $X$. These are equivalences in $\mfld^{\dd}$, as the composite $X \to Y \to X$ is isotopic to the identity, and likewise for $Y \to X \to Y$.   In particular, in $\mfld^{\dd}$, a point---which is equivalent to $[-1,1]^n$ by~\eqref{eqn. thickening X}---is equivalent to a disk of any dimension.
\end{remark}

\begin{remark}
\label{remark. frame bundle is morphism space}
Let $X = \RR^0$ be a point and $Y$ be any smooth compact manifold. Then the space of maps in $\mfld^{\dd}$ from $X$ to $Y$ is homotopy equivalent to the stable frame bundle of $Y$. To see this, note that if $X = [-1,1]^d$ and $\dim Y = d$, then the space of codimension zero smooth embeddings of $X$ into $Y$ is homotopy equivalent to the frame bundle of $Y$---i.e., the principle $GL(d)$-bundle associated to the tangent bundle. Further details are in the proof of Proposition~\ref{prop. hom from point is frame bundle}.
\end{remark}

The classifying map $X \to BGL \simeq BO$ of the stable tangent bundle survives the identification~\eqref{eqn. thickening X} up to homotopy. If a map $j: X \to Y$ is a smooth embedding of codimension zero, we further know $j$ respects the map to $BO$. We thus obtain a functor of $\infty$-categories
	\eqn\label{eqn. main functor}
	\mfld^{\dd}
	\to
	\Top_{/BO}^{\finite}
	\eqnd
to the $\infty$-category of spaces homotopy equivalent to a finite CW complex equipped with a map to $BO$. (One model of this functor is given in Construction~\ref{construction. main functor}.) 
In this work we give a new proof of the following fact:
	\begin{theorem}\label{theorem. main}
	\eqref{eqn. main functor} is an equivalence.
	\end{theorem}
This is one manifestation of the philosophy that questions of smooth embeddings become homotopy-theoretic upon thickening. Treated here is the case of smooth embeddings with trivialized normal bundles (i.e., thickened codimension-zero embeddings).

\begin{remark}[Precedents]
Theorem~\ref{theorem. main} is of no surprise to experts. Indeed, that~\eqref{eqn. main functor} is fully faithful can alternatively be proven using the $h$-principle. (Here is a sketch: The space of thickened embeddings might as well be a space of immersions with trivialized normal bundle, so invoke a version of the Smale-Hirsch theorem of immersions with trivialized normal bundles, increasing the dimensions of codomains.) 

At the level of isomorphism classes of both sides of~\eqref{eqn. main functor}---and in particular, considering the essential surjectivity of~\eqref{eqn. main functor}---there are many precedents. Barry Mazur's works~\cite{mazur-bulletin-1961-stable-equivalence-of-manifolds, mazur-1964-infinite-repetition} may have been the earliest. (Mazur proves stronger results about objects---up to diffeomorphisms, not up to isotopy equivalences---but weaker results in other respects, making no mention of the space of embeddings.) Mazur's results are referenced and enhanced by Waldhausen~\cite{waldhausen-1987-outline-how-manifolds-relate}, as an equivalence of the underlying $\infty$-groupoids in~\eqref{eqn. main functor} (rather than of the $\infty$-categories themselves). 
Wall~\cite{wall-1966-iv-thickenings} also investigated homotopy classes of maps to $BO$ in terms of what he calls ``thickenings'' of CW complexes. 
\end{remark}

\begin{remark}
The $\infty$-category of compact, smooth manifolds---possibly with corners and boundary, and allowing for all smooth maps---is equivalent to the $\infty$-category of spaces homotopy equivalent to finite CW complexes. By instead insisting on maps that are (thickened) codimension-zero embeddings (or, equivalently, smooth maps with trivialized normal bundles) we retain the data of the classifying map to $BO$, that is, the data of the stable tangent bundle.
\end{remark}

\begin{corollary}
Suppose two compact manifolds (possibly with corners) $X$ and $Y$ admit a {\em continuous} homotopy equivalence $X \to Y$ respecting the maps to $BO$ up to homotopy. Then there exists, possibly after thickening both $X$ and $Y$, a codimension zero smooth embedding from $X$ to $Y$, and from $Y$ to $X$, which are mutually inverse up to isotopy. 
\end{corollary}

\begin{corollary}
Suppose that $X$ is a compact smooth manifold. 
Then the space of maps from $X$ to $D^N$ in $\mfld^{\dd}$ is homotopy equivalent to the space of null-homotopies of the classifying map of the stable tangent bundle of $X$. 

In particular, if $X$ does not have stably trivial tangent bundle, the mapping space to $D^N$ is empty. Otherwise: The space of smooth embeddings into $D^N$ equipped with a trivialization of the normal bundle becomes, as $N$ goes to $\infty$, weakly homotopy equivalent to the space of continuous maps from $X$ to $O$.
\end{corollary}

The main geometric fact we use to prove Theorem~\ref{theorem. main} is the following:
	\begin{theorem}\label{theorem. pushouts exist}
	$\mfld^{\dd}$ has finite colimits. 
	\end{theorem}
While Theorem~\ref{theorem. pushouts exist} is obvious given Theorem~\ref{theorem. main},\footnote{We emphasize that in this note, Theorem~\ref{theorem. pushouts exist} is proven first.} at face value Theorem~\ref{theorem. pushouts exist} can perplex: We know full well that the collection of manifolds is not closed under gluing. Theorem~\ref{theorem. pushouts exist} states instead that the collection of thickened manifolds is closed under homotopy coherent gluing. Both thickening, and the $\infty$-categorical notion of (homotopy) colimit, are necessary.

We will exhibit the colimits explicitly. As a consequence, we will be able to interpret handle attachments as pushouts; this in turn allows us to show that a single object---the point---generates $\mfld^{\dd}$ under finite colimits (Lemma~\ref{lemma. point generates}). The constructions will also allow us to prove that the functor~\eqref{eqn. main functor} preserves finite colimits (Proposition~\ref{prop. main functor preserves finite colims}).

Now we may sketch the proof of Theorem~\ref{theorem. main}: Both the domain and codomain of~\eqref{eqn. main functor} are generated by the point, while the functor on the point is fully faithful by Remark~\ref{remark. frame bundle is morphism space}. More details are given in Section~\ref{section. proof of main theorem}.

Because $\mfld^{\dd}$ is a small $\infty$-category generated under finite colimits by a single object, Theorem~\ref{theorem. pushouts exist} makes formal that its Ind-completion is presentable. We may thus remove the finiteness in the codomain of Theorem~\ref{theorem. main}:

\begin{corollary}
\label{corollary. Ind mflds}
The induced functor
	\eqnn
	\Ind(\mfld^{\dd})
	\to
	\Top_{/BO}
	\eqnd
is an equivalence of $\infty$-categories.
\end{corollary}

\begin{example}[A model for the terminal object of $\Ind(\mfld^{\dd})$]\label{example. terminal object}
It is easily checked that $\mfld^{\dd}$ does not have a terminal object. For example, the space of maps from $X = \RR^0 \sim [-1,1]^N$ to any manifold $Y$ is homotopy equivalent to the stable frame bundle of $Y$ (Remark~\ref{remark. frame bundle is morphism space}). Any manifold whose stable frame bundle is contractible must be homotopy equivalent to $BO$---and $BO$ is not homotopy equivalent to the thickening of any compact manifold (with or without corners). After all, the cohomology of $BO$ is not bounded in degree.

Presentability guarantees a terminal object; the corollary tells us the terminal object has the homotopy type of $BO$. So there is an Ind-object of $\mfld^{\dd}$, homotopy equivalent to $BO$, which one may model as an increasing union of (thickened) manifolds, and which serves as a terminal object in $\Ind(\mfld^{\dd})$. 

The reader may be tempted to think that this presentation of $BO$ must be identical to the usual one: $BO$ is a colimit $\colim_{n,k} Gr_k(\RR^n)$ of Grassmannians, and the maps in the colimit diagram are smooth embeddings. However, these embeddings do not have trivial normal bundles, even stably. Indeed, when $k=1$, the embeddings $\RR P^n \to \RR P^{n+1}$ have normal bundles that pull back to the Mobius bundle along any $\RR P^1 \subset \RR P^n$, hence these normal bundles have non-trivial characteristic classes. 

What is true is that, in the usual diagram of Grassmannians, one may find a cofinal subdiagram wherein the normal bundle of each smooth embedding is trivializable, and by choosing trivializations, this cofinal subdiagram lifts to a diagram in $\Mfld^{\dd}$. To see this, recall that the tangent bundle of $Gr_k(\RR^n)$ is identified with
	\eqnn
	\hom(\gamma_{k,n},\gamma^{\perp}_{k,n}).
	\eqnd
Here, $\gamma_{k,n}$ is the tautological vector bundle on $Gr_{k}(\RR^n)$ whose fiber over $V \in Gr_k(\RR^n)$ is $V$ itself, and $\gamma_{k,n}^{\perp}$ has fibers given by the orthogonal complement of $V$ in $\RR^n$. For every $m  \geq j\geq 0$, the standard embedding $Gr_{k}(\RR^n) \to Gr_{k+j}(\RR^{n+m})$ thus has normal bundle
	\eqn\label{eqn. normal bundle to grassmannian embeddings}
	\hom(\underline{\RR^j},\gamma_{k,n}^{\perp})
	\oplus
	\hom(\gamma_{k,n},\underline{\RR^{m-j}})
	\oplus
	\hom(\underline{\RR^j},\underline{\RR^{m-j}})
	\cong
	(\gamma_{k,n}^{\perp})^{\oplus j}
	\oplus
	(\gamma_{k,n}^\vee)^{\oplus  m - j}
	\oplus
	\hom(\underline{\RR^j},\underline{\RR^{m-j}})
	\eqnd
where the underlines denote trivial vector bundles, and $\gamma_{k,n}^\vee$ denotes the $\RR$-linear dual vector bundle. Using the standard inner product on $\RR^n$, one has a short exact sequence of vector bundles on $Gr_k(\RR^n)$
	\eqnn
	0 \to \gamma_{k,n}^\perp \to \underline{\RR^n} \to \gamma_{k,n}^\vee \to 0
	\eqnd
and, because we are working continuously, this short exact sequence splits. In particular, so long as $m = 2j$, the normal bundle~\eqref{eqn. normal bundle to grassmannian embeddings} is trivializable. 

Thus, for any $0 \leq k \leq n$, the sequence of smooth manifold embeddings
	\eqnn
	Gr_{k}(\RR^n) 
	\to Gr_{k+1}(\RR^{n+2})
	\to Gr_{k+2}(\RR^{n+4})
	\to \ldots
	\eqnd
lifts to a sequence in $\Mfld^{\dd}$. By cofinality, the above sequence is a defining colimit for $BO$ in $\Top$.  This allows us to lift the sequence to $\Top_{/BO}$, where now the colimit is the terminal object of $\Top_{/BO}$ (see Remark~\ref{remark. slice categories create colimits}). By Corollary~\ref{corollary. Ind mflds}, the sequential diagram in $\Ind(\mfld^{\dd})$ has colimit given by the terminal object of $\Ind(\mfld^{\dd})$.
\end{example}

\subsection{Motivations}
\label{section. motivations}

As indicated above, Theorem~\ref{theorem. main} is a consequence of already available tools, but we had not seen the theorem in the literature. Not only do we record the result here, the proof in the present work is also distinct.
Indeed, a main motivation for this work was to give a proof of Theorem~\ref{theorem. pushouts exist} involving only straightforward constructions of manifolds, using no $h$-principles, and without invoking Theorem~\ref{theorem. main} (which one can now view as a consequence of thickened manifolds admitting finite colimits).
This was not for the sake of exercise. 

The author's main motivation was as follows.

It turns out that the $\infty$-category of Liouville sectors---a rich class of exact symplectic manifolds modeling cotangent bundles of singular spaces---also admits homotopy pushouts upon thickening. (This was mentioned, though not proven, in~\cite{last-stabilized}. Indeed, we learned the particular model for the pushout in Section~\ref{section. pushouts} from O. Lazarev.)

But there are fewer tools for dealing with embeddings of Liouville sectors, as the $h$-principle is not as robustly available, and there is no known analogue of Theorem~\ref{theorem. main}. It was with an eye toward understanding colimits of Liouville sectors that we endeavored to write this note. 

The central labor of the present work (Section~\ref{section. pushouts}) is the labor that carries over to the Liouville setting, and this is the only method we know for exhibiting colimits of Liouville sectors.

\begin{remark}
The  Ind-completion of the $\infty$-category of thickened Liouville sectors also has a terminal object, but we do not know its homotopy type at present. Compare with Example~\ref{example. terminal object}.
\end{remark}

\begin{remark}[We do not use any $h$-principles]
\label{remark. no h principles}
In a previous draft of this work, we came close to using an $h$-principle by invoking smooth approximation---i.e., that any continuous map between smooth manifolds may be continuously homotoped to be smooth---to prove that~\eqref{eqn. main functor}
is essentially surjective. 
As pointed out to us by Branko Juran, the use of smooth approximation is unnecessary. Indeed, in the same draft we had already proven that the functor preserves finite colimits and is fully faithful; and this is enough (see Section~\ref{section. proof of main theorem}).
\end{remark}

There is a second motivation to pursue an $h$-principle-free proof of Theorem~\ref{theorem. main}. As pointed out by an anonymous referee and by Branko Juran, it is natural to pursue analogues of Theorem~\ref{theorem. main} for PL manifolds and topological manifolds, where now the stabilized $\infty$-categories would be equivalent to $\Top^{\finite}_{/BPL}$ and $\Top^{\finite}_{/BTop}$, respectively. Indeed, if anything, all our finagling with smoothings becomes unnecessary. (Smoothing is the reason we have to deal with $P$ and $q$ as opposed to the more straightforward $P'$ and $q'$ in Section~\ref{section. pushouts}.) The main technical tools are still available via (i)  the theory of microbundles, and (ii) the existence of handle decompositions for PL manifolds (using stars) and for topological manifolds (upon thickening -- topological handle decompositions do not exist in all dimensions). 

Generalizing Theorem~\ref{theorem. main} to the PL and topological settings would produce rather clean formulations of not just smoothing at the object-level, but smoothing thickened embeddings by studying fibers of the functors $\Top^{\finite}_{/BO} \to \Top^{\finite}_{/BPL}$ and $\Top^{\finite}_{/BPL} \to \Top^{\finite}_{/BTop}$. (Example: One would have an obstruction-theoretic formulation of when a family of thickened codimension-zero topological embeddings lifts to a family of thickened codimension-zero smooth embeddings.) We note that the {\em object}-level question of when one can lift a thickened topological manifold to a thickened smooth manifold was already addressed completely in Milnor's original work on microbundles (see Theorem~5.13 of~\cite{milnor-microbundles}) with answers phrased in terms of topological $K$-theory with respect to $O$ and to $Top$.

We do not plan to pursue a full proof of the $PL$ and topological cases, nor a suitable exploration of these implications, at the moment.

\subsection{Convention: spaces, limits and colimits}
\label{section. conventions}
Throughout this note, we use the term ``mapping space'' to refer to a Kan complex of morphisms in a simplicially enriched category, or to a Kan complex of morphisms in an $\infty$-category.

Now-a-days it is sometimes common to use the term ``anima'' when a speaker chooses to remain agnostic about a model for a homotopy theory of spaces (so ``anima'' could refer to spaces, or to Kan complexes). However, in this work, we must use both Kan complexes and topological spaces (homotopy equivalent to CW complexes) explicitly. For example, our ``mapping spaces'' are typically Kan complexes, while the underlying space of a manifold is a topological space. Regardless, we do not encounter any topological space that cannot be recovered (up to homotopy equivalence) by its singular Kan complex, so this distinction will be immaterial for us when invoking $\infty$-categorical arguments---the anima aficionado may simply think of our models as anima.

We note that for $\cC$ an $\infty$-category and a diagram $I \to \cC$, there is only one notion of limit (and only one notion of colimit). This notion captures the intuition of what would traditionally be called a ``homotopy (co)limit.'' 

In contrast, for the category of topological spaces or of Kan complexes (and, more generally, any model category) there is utility in distinguishing the categorical notion of (co)limit from the $\infty$-categorical notion. So we may at times emphasize a ``point-set'' colimit of a diagram of topological spaces, for example in \eqref{eqn. diagram to P}. This means we take the colimit in the sense of the 1-category of topological spaces---i.e., the usual, non-homotopy-theoretic, notion of gluing topological spaces together. We will also use the term ``homotopy (co)limit'' in the context of topological spaces or Kan complexes for purposes of emphasis or disambiguation. See Proposition~\ref{prop. reduction to pullback of spaces} for an example of this paragraph and of the preceding paragraph.

\subsection{Acknowledgments} 
The author was supported by an NSF CAREER grant DMS-2044557, an Alfred P. Sloan Research Fellowship, a Texas State University Presidential Seminar Award and Valero Award, and a Research Membership at MSRI/SLMATH. We thank 
Diarmuid Crowley,
John Francis,
S\o ren Galatius,
John Klein,
Sander Kupers, and
Jacob Lurie
for communications educating me on past works. We also thank Branko Juran for his interest and insight (see Remark~\ref{remark. no h principles}).
Finally, we thank the anonymous referees for their helpful comments.

\section{From smooth manifolds to spaces over \texorpdfstring{$BO$}{BO}}\label{section. from mfld to TopBO}
Here we define $\mfld^{\dd}$ and give a construction of~\eqref{eqn. main functor}.

\subsection{Preliminaries on manifolds with corners}

\begin{remark}[Conventions for manifolds with corners]
There are various definitions of manifolds with corners in the literature, especially when one begins to endow manifolds with stratifying data. We supply no such data and we use the minimal definition of smooth manifolds with corners. We refer to the first paragraphs of Section~2 of~\cite{joyce-manifolds-with-corners}; we also note that our notion of smoothness is called {\em weakly smooth} there.

In particular, an $n$-manifold with corners $X$ in our work is a paracompact\footnote{We note  that all the manifolds we consider in our categories will be compact; but the use of local arguments makes it convenient to define manifolds with corners in non-compact settings. Our manifolds will also have components of equal dimension later on; this is irrelevant for the present, general discussion.} Hausdorff topological space equipped with a maximal atlas of smoothly compatible local charts from open subsets of the octant $(\RR_{\geq 0})^n$. Note that for any $x \in X$, whether $x$ is the image of a codimension $k$ face of the octant is invariant under change of charts.

As usual, a map between Euclidean octants is called smooth if it arises as the restriction of a smooth map defined on a small neighborhood of the octants. Through charts, this defines a notion of smoothness for maps between manifolds with corners in the usual way. Note in particular that smooth maps between manifolds with corners, in this work, need not respect boundary/corner strata---for example, a compact disk may smoothly map to Euclidean space.
\end{remark}

\begin{defn}
A smooth embedding $j: X \to Y$ is called an {\em isotopy equivalence} if there exists a smooth embedding $h: Y \to X$ together with smooth isotopies $jh \sim \id_Y$ and $hj \sim \id_X$. (Just as with $h$ and $j$, the smooth isotopies need not respect strata in any way.)
\end{defn}

\begin{remark}[One can replace corners with boundaries]\label{remark. isotopy equivalence removes corners}
Consider a compact smooth manifold $X$ with non-empty corner and boundary strata. There exists a smooth manifold with boundary (but no corners) $X'$ and an isotopy equivalence $X' \into X$. 
Details for one such construction may be found in Douady and H\'erault's appendix to Borel and Serre's compactification paper~\cite{borel-serre-compactification}.

We do caution that, it seems to this writer that $X'$ is not unique up to diffeomorphism---indeed, two such $X'$ must have smoothly cobordant boundary, but it is not clear that even the boundaries of $X'$ need be diffeomorphic (though Douady and H\'erault's methods will produce {\em homeomorphic} boundaries). We have not come up with an example showing that the $X'$ need not be diffeomorphic, though we strongly suspect there exist such examples. It is clear that $X'$ is unique up to isotopy equivalence (as all $X'$ are isotopy equivalent to $X$); this suffices for us.
\end{remark}

\begin{remark}
One may, in light of Remark~\ref{remark. isotopy equivalence removes corners}, be tempted to consider only smooth manifolds with boundary (and no corners) but this would complicate the coherence of the thickening process that follows. So we allow for corners.
\end{remark}

\subsection{Defining the \texorpdfstring{$\infty$}{infinity}-category of thickened manifolds}

\begin{notation}[$\MMfld_d$]\label{notation. mfld d}
Let $\MMfld_d$ denote the Kan-complex enriched category whose objects are compact, smooth manifolds of dimension $d$, and whose morphism spaces are (the singular complex of) the spaces of smooth embeddings. (Note we do not demand that the embeddings respect corner strata in any way.) 
\end{notation}

\begin{remark}
One model for the Kan complex of morphisms from $X$ to $Y$ is as follows. Let $\Delta^k_e$ be the subspace of $\RR^{k+1}$ whose coordinates sum to 1. A $k$-simplex of $\hom_{\MMfld_d}(X,Y)$ is the data of a smooth embedding $X \times \Delta^k_e \to Y \times \Delta^k_e$ respecting the projections to $\Delta^k_e$. This was also used for example in~\cite[Convention~4.1.5]{aft-1}. 
\end{remark}

\begin{remark}
The natural notion of equivalence (i.e., a map that admits an inverse up to homotopy) in $\MMfld_d$ coincides with the notion of isotopy equivalence.
\end{remark}

We have a functor of simplicially enriched categories
	\eqn\label{eqn. thickening functor}
	\MMfld_d \to \MMfld_{d+1},
	\qquad
	X \mapsto X \times [-1,1]
	\eqnd
taking a manifold to a direct product with the compact unit interval. 

\begin{notation}[$\MMfld^{\dd}$]
\label{notation. thickened manifold category}
We define
	\eqn\label{eqn. mmfld increasing union}
	\MMfld^{\dd} := \colim\left(
	\MMfld_0 \xrightarrow{-\times [-1,1]} \MMfld_1 \xrightarrow{-\times [-1,1]} \ldots
	\right),
	\eqnd
concretely modeled as an increasing union of simplicially enriched categories.\footnote{This is also a homotopy colimit in the model category of simplicially enriched categories.} We call $\MMfld^{\dd}$ the (simplicially enriched) category of thickened compact manifolds. 
\end{notation}

Informally, $\MMfld^{\dd}$ is a category where any compact manifold $X$ is identified with $X \times [-1,1]^N$ for any $N\geq 0$, and where a morphism from $X$ to $Y$ is a codimension zero embedding $X\times[-1,1]^N \to Y \times [-1,1]^{N'}$ for some $N,N'$.

\begin{notation}[$\mfld^{\dd}$]\label{notation. mfld thickened}
We define the $\infty$-category
	\eqnn
	\mfld^{\dd} := N(\MMfld^{\dd})
	\eqnd
to be the homotopy coherent nerve of the simplicially enriched category $\MMfld^{\dd}$.
\end{notation}

\begin{remark}
\label{remark. hom spaces in mfld and MMfld}
We recall that the homotopy coherent nerve $N$ commutes with filtered colimits, so one could equivalently define $\mfld^{\dd}$ as the increasing union of $\infty$-categories 
	\eqnn
	N(\MMfld_0) \to N(\MMfld^{1}) \to \ldots
	\eqnd
obtained by applying $N$ to~\eqref{eqn. mmfld increasing union}.

We also note that the mapping spaces in the nerve are homotopy equivalent to the mapping spaces of the simplicially enriched category, so long as the simplicially enriched category has morphisms given by Kan complexes (this hypothesis is met for $\MMfld_d$ for all $d$ and for $\MMfld^{\dd}$)---see for Example Theorems~1.1.5.13 and 2.2.0.1 of~\cite{lurie-htt}. As a result, for any two objects $X,Y \in \MMfld^{\dd}$, we have a homotopy equivalence of Kan complexes
	\eqn\label{eqn. hom mmfld vs hom mfld}
	\hom_{\mfld^{\dd}}(X,Y)
	\simeq
	\hom_{\MMfld^{\dd}}(X,Y).
	\eqnd
\end{remark}

Any manifold $X$ of dimension $d$ is the base of the frame bundle $\Fr_d(X)$, which is a principle $GL(d)$-bundle. By thickening, we obtain a space $\Fr(X)$ with a free $GL$-action, where $GL$ is the infinite general linear group, which we call the stabilized frame bundle of $X$. Let $\TTop^{GL}$ denote the simplicially enriched category of topological spaces (homotopy equivalent to CW complexes) equipped with a continuous $GL$-action. Then the stable frame bundle construction gives rise to a functor of $\infty$-categories
	\eqn\label{eqn. Fr}
	\Fr: \Mfld^{\dd} \to N(\TTop^{GL})
	\eqnd
to the $\infty$-category of spaces with $GL$-action (obtained as the homotopy coherent nerve of $\TTop^{GL}$).

\begin{remark}
\label{remark. model for Fr}
A concrete model of $\Fr$ may be given as follows. Note the commutative (up to natural transformation) diagram of simplicially enriched categories
	\eqn\label{eqn. thickening frame bundles}
	\xymatrix{
	\MMfld_d \ar[r]^{\Fr_d} \ar[d]_{-\times [-1,1]}
		& \TTop^{GL(d)} \ar[d]^{-\times_{GL(d)} GL(d+1)} \\
	\MMfld_{d+1} \ar[r]^{\Fr_{d+1}}
		& \TTop^{GL(d+1)}
	}
	\eqnd
where $\Fr_{d} \times_{GL(d)} GL(d+1)$ is the principal $GL(d+1)$-bundle associated to $\Fr_d$ via the group homomorphism $GL(d) \to GL(d+1)$. Because the natural transformation maps are---for all objects $X \in \MMfld_d$---equivalences in $N(\TTop^{GL(d+1)})$, the homotopy-coherent nerve $N$ renders~\eqref{eqn. thickening frame bundles} to a diagram $\Delta^1 \times \Delta^1 \to \inftycat$ in the $\infty$-category of $\infty$-categories. Noting the natural map 
	\eqnn
	\colim(\TTop^{GL(0)} \to \TTop^{GL(1)} \to \TTop^{GL(2)} \to \ldots )
	\to \TTop^{GL},
	\eqnd
and noting that $N$ commutes with filtered colimits, the colimit (indexed by $d$) of~\eqref{eqn. thickening frame bundles} induces~\eqref{eqn. Fr}.
\end{remark}

\subsection{The functor to spaces over \texorpdfstring{$BO$}{BO}}

\begin{remark}
It is classical that the maps $O(d) \to GL(d)$ are homotopy equivalences of group-like $E_1$-spaces, and that the induced map $O \to GL$ is an equivalence of group-like $E_1$ spaces. (In fact, though we will not need this, this map may be promoted to an equivalence of infinite loop spaces.) As a result, the induced functor (restricting a $GL$-action to an $O$-action)
	\eqn\label{eqn. GL spaces to O spaces}
	N(\TTop^{GL}) \to N(\TTop^{O})
	\eqnd
is an equivalence of $\infty$-categories.

Let $BO$ denote the Kan complex modeling the classifying space of the simplicial group $O$. (The interested reader may find a concrete model for $BO$ and $EO$ in~\cite[page 87]{may-simplicial-objects-1968}, where $BO$ is denoted by $\overline{W}(H)$ with $H=O$.) It is classical that one has a functor of $\infty$-categories
	\eqn\label{eqn. functor from O-spaces to over BO}
	N(\TTop^O) \to \Top_{/BO}
	\eqnd
from the nerve of the simplicially enriched category $\TTop^O$ to the $\infty$-category $\Top_{/BO}$ of spaces (Kan complexes) equipped with a map to $BO$. Informally, the map~\eqref{eqn. functor from O-spaces to over BO} takes an $O$-space, equipped with the $O$-equivariant map to a point, to its homotopy quotient, equipped with the induced map to $pt/O \simeq BO$.

Moreover, \eqref{eqn. functor from O-spaces to over BO} is an equivalence of $\infty$-categories. A proof using simplicially enriched model categories (where equivalences between $O$-spaces are the ``coarse'', also known as ``naive,'' equivalences) goes back to~\cite{dror-dwyer-kan-1980}.
\end{remark}

\begin{construction}\label{construction. main functor}
Composing~\eqref{eqn. Fr} and~\eqref{eqn. GL spaces to O spaces} and~\eqref{eqn. functor from O-spaces to over BO}, we obtain a functor of $\infty$-categories
	\eqnn
	\mfld^{\dd}
	\to
	\Top_{/BO}.
	\eqnd
In fact, it is easy to note that (because all of our manifolds are compact) the above functor factors through the $\infty$-category of (spaces homotopy equivalent to) finite CW complexes equipped with maps to $BO$:
	\eqnn
	\mfld^{\dd}
	\to
	(\Top^{\finite})_{/BO}.
	\eqnd
This is the map~\eqref{eqn. main functor}.
\end{construction}

\begin{example}
Given a manifold $X$ of dimension $d$, its image in $(\Top^{\finite})_{/BO}$ is computed as follows.
First, $\Fr(X)$ is a space with a free $O = O(\infty)$ action. The image of this space under~\eqref{eqn. functor from O-spaces to over BO} is the homotopy quotient of $\Fr(X)$ by the $O$-action, equipped with its natural map to $BO$. In particular,~\eqref{eqn. functor from O-spaces to over BO} sends $X$ to a space homotopy equivalent to $X$ equipped with a map to $BO$. We know this map is homotopic to the map classifying the stable frame bundle, as by construction it fits into a homotopy pullback square of spaces
	\eqnn
	\xymatrix{
	\Fr(X) \ar[r] \ar[d] & \ast \simeq EO \ar[d] \\
	X \ar[r] & BO
	}
	\eqnd
where both vertical arrows are quotients by $O$.
\end{example}

\begin{remark}[Alternative construction of~\eqref{eqn. main functor}]
\label{remark. alternative construction}
As shown in the proof of Proposition~\ref{prop. hom from point is frame bundle}, one has a natural identification of $\hom_{\mfld^{\dd}}(pt,X)$ with the stable frame bundle $\Fr(X)$ -- as an $O$-space over $X$. Thus, as pointed out to us by an anonymous referee, the assignment $X \mapsto \hom_{\mfld^{\dd}}(pt,X)$ gives an alternative characterization of~\eqref{eqn. main functor}.
\end{remark}

\section{Pushouts of thickened manifolds}
\label{section. pushouts}

\begin{setup}\label{setup. WXY}
Throughout, we fix three compact manifolds $W, X, Y$. 
Stabilizing as necessary, we assume $\dim X = \dim W = \dim Y$. 
By Remark~\ref{remark. isotopy equivalence removes corners}, we may further assume that $W,X,Y$ have no corners (but they may have boundary). We equip them with smooth, codimension zero embeddings
	\eqn\label{eqn. iX iY}
	X \xleftarrow{i_X} W \xrightarrow{i_Y} Y.
	\eqnd
\end{setup}

\begin{example}\label{example. W has no boundary}
In the special case that $W$ is a manifold with no boundary and no corners, we see that $i_X$ and $i_Y$ are diffeomorphisms onto certain connected components of $X$ and $Y$. 

(To see this, note that the interior of $W$---meaning the subspace of $W$ admitting charts from $\RR^n$, or equivalently, the complement of the boundary and corner strata---is an open subset of $W$. If $W$ has no boundary or corners, then $W$ is equal to its interior. On the other hand, any smooth embedding from a manifold with no boundary and no corners to a manifold of the same dimension is an open mapping. Thus when $W$ is compact the image of $i_X$ is both open and closed.)
\end{example}

\subsection{The model \texorpdfstring{$P$}{P} of a pushout}
\label{section. construction of P}
Given~\eqref{eqn. iX iY}, consider the topological space
	\eqn\label{eqn. P before smoothing}
	P' := P'(i_X,i_Y) :=
	[-2,-1] \times X \bigcup_{i_X \times \id_{\{-1\}}} [-1,1] \times W \bigcup_{i_Y \times \id_{\{1\}}} [1,2] \times Y
	\eqnd
(which one can recognize as one model for a homotopy pushout of spaces). See Figure~\ref{figure. P and P'}.

\begin{remark}
When $W$ has boundary, $P'$ is not canonically a smooth manifold with corners. This is because one must choose a smooth atlas near $\{-1,1\} \times \del W$. On the other hand, if we remove the locus $\{-1,1\} \times \del W$, $P'$ has a canonical smooth structure. This is a standard issue. The smooth structure on $P'$ is only as undefined as a smooth structure on the (closed) complement of a(n open) quadrant in $\RR^2$. Indeed, an outward normal to $\del W$ and normal vectors to $\{-1,1\}$ inside $[-2,2]$ reduce the problem to $(\RR^2\setminus\text{Quadrant}) \subset \RR^2$ by considering the embedding
	\eqnn
	(\RR^2\setminus\text{Quadrant}) \times \del W \into P'.
	\eqnd
\end{remark}

In the following remarks, we present two ways to treat $P'$ as, or alter $P'$ to be, a smooth object. We are careful here because, to avoid notational clutter, it will be very convenient to be able to treat the obvious embeddings
	\eqn\label{eqn. the maps we want to be smooth}
	[-2,-1] \times X \into P',
	\qquad
	[-2,2] \times W \into P',
	\qquad
	[1,2] \times Y \into P'
	\eqnd
as all smooth. (Note that the embedding $[-2,2] \times W \into P'$ is set-theoretically defined piece-wise, equaling $i_X$ and $i_Y$ along the subintervals $[-2,-1]$ and $[1,2]$, respectively.)

\subsection{A smoothing of \texorpdfstring{$P'$}{P'}}
\label{section. smoothing P of P'}

\begin{choice}
Choose a small positive real number $\epsilon$. (All that matters is that $\epsilon$ is less than 1, and this 1 is only relevant because we play with intervals of length 1 below.)
\end{choice}

Note we may also present $P'$ as the point-set colimit of the diagram of topological spaces
	\eqn\label{eqn. diagram to P}
	\xymatrix{
			&&W \times [1,1+\epsilon]
			\ar[d]^{\id_W \times \iota}
			\ar[rr]^{i_Y \times \iota} 
				&& Y \times [1,2] \\
	W \times [-1-\epsilon,-1]  
		\ar[rr]^{\id_W \times \iota}
		\ar[d]^{i_X \times \iota}			
		&&W \times [-1-\epsilon,1+\epsilon] \\
	X \times [-2,-1]
	}
	\eqnd
where the $\iota$ maps are inclusions of subintervals.

\begin{notation}[$u$]
We will denote an element of (all of) the intervals in~\eqref{eqn. diagram to P} by $u$.
\end{notation}

\begin{example}\label{example. P' is smooth if W is closed}
In the setting of Example~\ref{example. W has no boundary}, $P'$ inherits an obvious smooth atlas and is hence a smooth manifold (possibly with boundary and corners) with atlas induced by those of $W, X, Y$ and $[0,1]$. If $W,X,Y$ are furthermore connected, then $P'$ is diffeomorphic to $X \times [-2,2] \cong X \times [0,1]$.
\end{example}

We make a smooth model $P$ of $P'$ as follows.

\begin{choice}[A coordinatized neighborhood of $\del W$ in $P'$]
\label{choice. nbhd of del W}
As per Set-up~\ref{setup. WXY}, we have assumed that $W,X,Y$ have no corners (but they may have boundary). This means $W$ admits a vector field strictly outward-pointing along its boundary. So choose one. Flowing along this vector field for time $-2\epsilon$ (which is possible because $W$ is compact) 
realizes a self-embedding 
	\eqnn
	c: W \to W
	\eqnd
whose image is disjoint from $\del W$. 
This induces a diffeomorphism
	\eqn\label{eqn. diff coordinatizing v}
	W \bigcup_{\del W \times (-\infty,0]_v} \left(\del W \times (-\infty,\epsilon)_v\right)
	\xrightarrow{\cong}
	W \setminus \del W.
	\eqnd
The subscript $v$ indicates the variable we will use for elements of the intervals in~\eqref{eqn. diff coordinatizing v}; note we have parametrized so that the $v=-\epsilon$ locus in the domain is sent to $c(\del W)$ in the codomain.
\end{choice}

\begin{remark}\label{remark. identifying neighborhoods in P'}
We may consider the direct product of~\eqref{eqn. diff coordinatizing v} with $\iota$ to think of $(\del W \times (-\infty,\epsilon)_v) \times [-2,2]_u \cong \del W \times [-2,2]_u \times (-\infty,\epsilon)_v$ as a subset of $P'$. 
Doing so, inside $P'$ there is a neighborhood of $\del W \times \{-1\}_u \times \{0\}_v$ homeomorphic to
	\eqnn
	\del W \times \left(
	(-1-\epsilon,-1+\epsilon)_u \times (-\infty,0]_v \bigcup (-1-\epsilon,-1]_u \times (-\infty,\epsilon)_v
	\right).
	\eqnd
Likewise we may think of $\del W \times \{1\}_u \times (-\infty,\epsilon)_v$ as a subset of $P'$ and we obtain an identification 
	\eqnn
	\nbhd\left(\del W \times \{1\}_u \times \{0\}_v\right)
	\cong
	\del W \times \left(
	(1-\epsilon,1+\epsilon)_u \times (-\infty,0]_v \bigcup [1,1+\epsilon)_u \times (-\infty,\epsilon)_v
	\right).
	\eqnd
\end{remark}

\begin{figure}[ht]
    \begin{equation}\nonumber
			\xy
			\xyimport(8,8)(0,0){\includegraphics[width=6in]{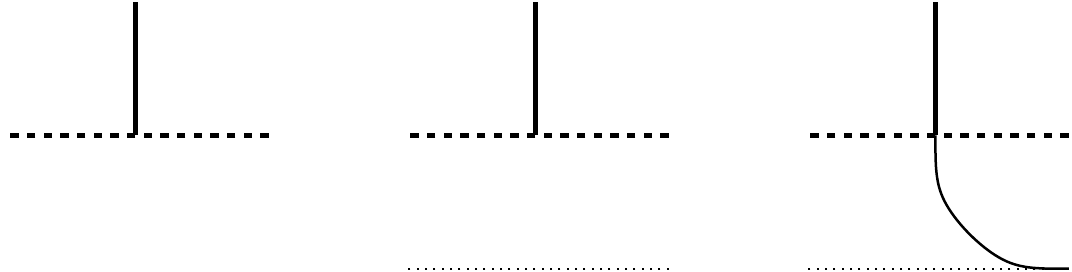}}
				,(1,-1)*{(i)};
				,(4,-1)*{(ii)};
				,(7,-1)*{(iii)};
				,(4.4,8)*{u=-1};
				,(3.1,4.5)*{v=0};
				,(3.1,0.5)*{v=-\epsilon};
			\endxy
    \end{equation}
    \caption{
    Illustrations of Choice~\ref{choice. smoothing corner}. In all illustrations is a neighborhood of $\del W \times \{-1\}_u \times \{0\}_v$ inside $P'$, projected onto the square $(-1-\epsilon,-1+\epsilon)_u \times [-\epsilon,\epsilon)_v$. As in Figure~\ref{figure. P and P'}, the $u$ direction is drawn horizontally. The $v$ direction is drawn vertically. In $(i)$, $\del W \times \{-1\}_u \times \{0\}_v$ is the fiber over the point where the solid line intersects the dashed line. The fiber above the dashed line is the locus $\del W \times (-1-\epsilon,-1+\epsilon)_u \times \{0\}_v$. The region below the dashed locus represents the region $\del W \times (-1-\epsilon,-1+\epsilon)_u \times [-\epsilon,0)_v$.
    In $(ii)$  is now a thin dotted line, indicating the locus $\del W \times (-1-\epsilon,-1+\epsilon)_u \times \{ -\epsilon\}_v$.
   The new curve in $(iii)$  indicates $\gamma$. The region consisting of points to the left of or below $\gamma$ is $R_X$. Note that $R_X$ intersects some portion of the thick dashed line, but does not contain all of the thick dashed line.
	}
    \label{figure. gamma}
\end{figure}

\begin{choice}[Smoothing a concave corner using $\gamma,R_X,R_Y$]
\label{choice. smoothing corner}
Now choose a set
	\eqnn
	\gamma \subset A : = (-1-\epsilon,-1+\epsilon)_u \times (-2\epsilon,\epsilon)_v
	\eqnd
such that
\enum[(i)]
	\item $\gamma$ is a smooth, connected 1-dimensional submanifold. 
	\item For some small ball $B$ containing the point $\{-1\}_u \times \{-\epsilon\}_v \in A$, we have 
		\eqnn
		\gamma \bigcap (A \setminus B) = 
		\left(
		[-1,\infty)_u \times \{-\epsilon\}_v \bigcup \{-1\}_u \times [-\epsilon,\infty)_v 
		\right)
		\bigcap (A \setminus B).
		\eqnd
	That is, outside of a tiny neighborhood of $\{-1\}_u \times \{-\epsilon\}_v$, $\gamma$ is equal to two positive rays.	
	\item $\gamma$ is contained in the region where $u\geq -1$ and $v \geq -\epsilon$.
\enumd
See Figure~\ref{figure. gamma}.
By (a version of) the Jordan Curve Theorem, $\gamma$ divides $A$ into two connected regions. We let 
	\eqnn
	R_X\subset A
	\eqnd
denote the region containing $\{-1\}_u \times \{-\epsilon\}_v$ and with boundary given by $\gamma$, so $R_X$ is a smooth 2-manifold with smooth boundary. We informally think of $R_X$ as smoothing the region $\{u \leq -1\} \bigcup \{v \leq -\epsilon\}$ (or as smoothing $\{u \leq -1\} \bigcup \{v \leq 0\}$ by removing a chunk of it).

Likewise, there is a smooth manifold with boundary 
	\eqnn
	R_Y \subset (1-\epsilon,1+\epsilon)_u \times (-2\epsilon,\epsilon)_v
	\eqnd
smoothing the locus $\{u \geq 1 \} \bigcup \{v \leq -\epsilon\}.$
\end{choice}

\clearpage

\begin{figure}[ht]
    \begin{equation}\nonumber
			\xy
			\xyimport(8,8)(0,0){\includegraphics[width=6in]{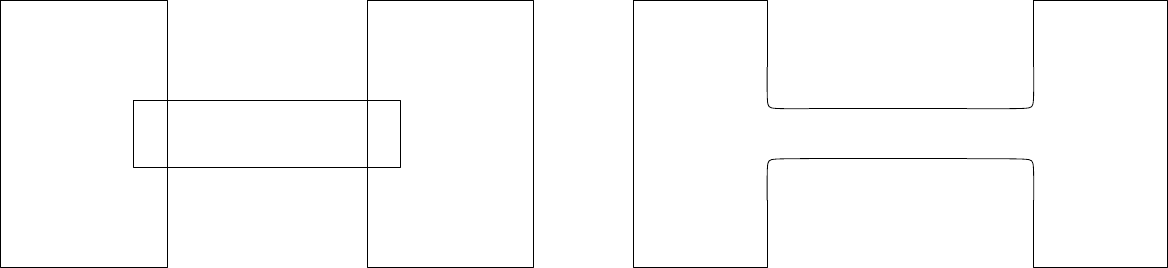}}
				,(0.5,-1)*{X \times [-2,-1]};
				,(1.85,3.9)*{W \times[-1-\epsilon,1+\epsilon]};
				,(3.1,-1)*{Y \times [1,2]};
				,(4.8,-1)*{X \times [-2,-1]};
				,(7.55,-1)*{Y \times [1,2]};
			\endxy
    \end{equation}
	\caption{
    A cartoon of $P'$ on the left (with regions indicated), and of $P$ on the right.  The intervals runs horizontally, while the vertical directions are meant to convey movement in the $W, X, Y$ directions. We have drawn not only $W \times [-1,1]$, but the whole $W \times [-1-\epsilon,1+\epsilon]$ to emphasize that $W$ admits a collar inside $X \times [-2,-1]$ and inside $Y \times [1,2]$. The drawing of $P$ (note that the cornered are rounded) is meant to indicate that $P$ is obtained by removing the not-obviously-smooth regions where the rectangles in the $P'$ drawing intersect.
	}
    \label{figure. P and P'}
\end{figure}

\clearpage

\begin{construction}[$P$]\label{construction. P}
We let $P \subset P'$ denote the smooth manifold with boundary (and no corners)
	\eqnn
	P:= 
	\left( X \times [-2,-1]_u \right)
	\bigcup 
	\left( \del W \times R_X  \right)
	\bigcup 
	\left( c(W) \times [-1-\epsilon,1+\epsilon]_u  \right)
	\bigcup 
	\left( \del W \times R_Y  \right)
	\bigcup 
	\left( Y \times [1,2]_u \right).
	\eqnd
See Figure~\ref{figure. P and P'}. Note that we are using Remark~\ref{remark. identifying neighborhoods in P'} to identify $\del W \times R_X$ and $\del W \times R_Y$ as subsets of $P'$. We think of $c(W) \times [-1-\epsilon,1+\epsilon]$ as a subset of $W \times [-1-\epsilon,1+\epsilon]$ in the defining diagram of $P$~\eqref{eqn. diagram to P}.
\end{construction}

\begin{notation}[$\alpha$]
By construction, $P$ is a subset of $P'$. The inclusion defines a smooth embedding 
	\eqn\label{eqn. alpha from P to P'}
	\alpha: P \into P'
	\eqnd
where the smoothness of $a$ makes sense because $a$ avoids the locus $\{-1,1\} \times \del W$.
\end{notation}

\begin{choice}[$\beta$]
Choose a smooth function $\phi$ compactly supported inside $\del W \times (R_X \coprod R_Y) \subset P'$, and consider the vector field $-\phi \del_v$ defined on all of $P'$. For a well-chosen $\phi$, we can arrange for the time-1 flow to have image completely contained in $P$. 
This isotopy is clearly smooth in the coordinate system invoked in Remark~\ref{remark. identifying neighborhoods in P'}. We will denote by 
	\eqnn
	\beta: P' \into P
	\eqnd
the end result (i.e., the time-1 map) of a choice of such an isotopy.
\end{choice}

\begin{remark}\label{remark. beta respects u coordinate}
Note that, by construction, $\beta$ respects the projection to $[-2,2]_u$. 
\end{remark}

\begin{remark}
Note that $\beta\alpha$ is smoothly isotopic to $\id_P$. Also, even though $P'$ is not a smooth manifold a prior, let us define a map out of $P'$ to be smooth if and only if the pre-composition with $\beta^{-1}: \beta(P') \to P'$ is smooth. Then one also finds that $\alpha\beta$ is smoothly isotopic to $\id_{P'}$. Thus, for any test object $Z$---that is, any object $Z$ of  $\MMfld_{\dim P}$---we find a  homotopy equivalence
	\eqn\label{eqn. emb P is emb'P'}
	\alpha^*: \sing \emb'(P',Z)  \to \sing\emb(P,Z)
	\eqnd
where $\emb'$ is the space of smooth embeddings out of $P'$ (in the sense defined in this remark). By definition, the maps in~\eqref{eqn. the maps we want to be smooth} are smooth embeddings; in this way, we conclude that $P$ is a pushout in the $\infty$-category $\mfld^{\dd}$ if $\sing\emb'(P',-)$---upon passing to $P' \times [0,1]^k$ as $k \to \infty$---models a functor out of $\mfld^{\dd}$ corepresented by a pushout.
\end{remark}

\begin{remark}
The analogues of~\eqref{eqn. the maps we want to be smooth}---obtained by post-composing~\eqref{eqn. the maps we want to be smooth} with $\beta$---
	\eqn\label{eqn. maps from X W Y into P}
	[-2,-1] \times X \into P,
	\qquad
	[-2,2] \times W \into P,
	\qquad
	[1,2] \times Y \into P,
	\eqnd
are all smooth, codimension-zero embeddings.
\end{remark}

The following is easy to verify; ``smoothing then thickening'' is isotopy equivalent to ``thickening then smoothing.'' We omit the proof.

\begin{prop}\label{prop. thickening P}
There are natural isotopy equivalences
	\eqnn
	P(i_X \times \id_{[0,1]} ,i_Y\times \id_{[0,1]} ) \to P(i_X,i_Y)\times [0,1].
	\eqnd
\end{prop}

\begin{remark}
Proposition~\ref{prop. thickening P} states that the construction of $P$ is compatible with thickening. So the restriction map along~\eqref{eqn. maps from X W Y into P} fits into a homotopy-commutative diagram of Kan complexes
	\eqnn
	\xymatrix{
	\hom_{\MMfld_{\dim P}}(P,-) \ar[r] \ar[d] & \hom_{\MMfld_{\dim P}}([-2,-1] \times X, -) \ar[d]  \\
	\hom_{\MMfld_{\dim P+1}}(P \times [0,1],- \times [0,1]) \ar[r] & \hom_{\MMfld_{\dim P+1}}([-2,-1] \times X \times [0,1], - \times [0,1]) 
	}
	\eqnd
and hence passes to the $\infty$-category of thickened manifolds
	\eqnn
	\hom_{\mfld^{\dd}}(P,-) \to \hom_{\mfld^{\dd}}([-2,-1] \times X,-)
	\eqnd
(and likewise for $[1,2] \times Y$ and $[-2,2] \times W$).
\end{remark}

\begin{remark}[Another approach using diffeological spaces]
\label{remark. P as a diffeological space}
We may also expand our notion of smooth manifold-with-corners, and treat $P'$ as a diffeological space---i.e., a space where we know what we mean by a smooth map into $P'$. Indeed, there is a natural notion of smooth map into $P'$. Consider the map $\text{Quadrant} \times \del W \into P'$ induced by the normal vectors from the previous paragraph. Then we say a map $f: A \to P'$ is smooth if (i) when pulled back along the inclusion $\text{Quadrant} \times \del W \into P'$, the composition with $\text{Quadrant} \times \del W \into \RR^2 \times \del W$ is smooth, and (ii) $f$ is smooth when pulled back along the complement of $\{-1,1\} \times W$. It is easy to see this is independent of choice of normal coordinates. 

For this choice of smooth structure on $P'$, one can arrange for our smoothing $P$ from Construction~\ref{construction. P} to sit into a diagram $P \into P' \into P \into P'$ where every composition is isotopic to the identity morphism---that is, $P \to P$ is isotopic to the identity of $P$ through smooth embeddings, and $P' \to P'$ is isotopic to the identity through diffeologically smooth embeddings. 
Thus, we may simply treat $P'$ as a smooth manifold isotopy equivalent to an honest smooth manifold $P$ possibly with corners. We do not take this approach here, at the expense of cluttering some formulas with $a$ and $b$.
\end{remark}

\subsection{Reduction to computing a homotopy pullback of spaces}

\begin{notation}[$I_s$]
\label{notation. I_s}
For later notational clarity, we let 
	\eqnn
	I_s := [0,1]
	\eqnd denote the interval by which we thicken manifolds. The subscript $s$ is to emphasize that we will denote elements of $I_s$ by $s$.
\end{notation}

\begin{choice}[$\sigma,h_t$]
\label{choice. h_t}
Here, we endow every interval in $\RR$ with the standard orientation inherited from $\RR$.
We fix the following:
\enum[(a)]
	\item An orientation-preserving diffeomorphism
		\eqnn
		\sigma: [-2,2] \to I_s.
		\eqnd
	\item A smooth isotopy of embeddings
		\eqnn
		\{h_t: I_s \to [-2,2]\}_{t \in [-2,2]}
		\eqnd
	satisfying the following properties:
		\enum[(i)]
		\item $h_t = h_{-1}$ for all $t \leq -1$.
		\item $h_t = h_{1}$ for all $t \geq 1$. (These first two conditions amount to a collaring condition on the isotopy.)
		\item $h_{-1}$ is an orientation-preserving diffeomorphism onto the image $[-2,-1]$.
		\item $h_{1}$ is an orientation-preserving diffeomorphism onto the image $[1,2]$.
		\item $h_t$ is strictly increasing for $-1 < t < 1$, meaning that for every $s \in I_s$, the derivative ${\frac{\del}{\del t}}h_t(s)$ is positive.
		\enumd
\enumd
\end{choice}

Using our choice of isotopies $h_t$, we thus have a diagram $\Delta^1 \times \Delta^1 \to N(\MMfld_{\dim P})$ we informally depict as follows:
	\eqn\label{eqn. first pushout}
	\xymatrix{
	I_s \times W \ar[r]^{h_{1} \times i_Y} \ar[d]^{h_{-1} \times i_X}
		& [1,2] \times Y \ar[d] \\
	[-2,-1] \times X \ar[r] 
		& P
	}.
	\eqnd
Here,  the $i_X,i_Y$ are as in~\eqref{eqn. iX iY}, and unlabeled maps are those from~\eqref{eqn. maps from X W Y into P}. The diagram above commutes up to the isotopy parametrized by $t \in [-2,2]$ given as follows:
	\eqnn
	\{ I_s \times W \xrightarrow{h_t \times \id_W} [-2,2] \times W \xrightarrow{\eqref{eqn. maps from X W Y into P}} P\}_{t \in [-2,2]}.
	\eqnd
Thus, for any test object $\bullet \in \MMfld^{\dd}$, by applying the functor $\hom_{\MMfld^{\dd}}(-,\bullet)$ (and invoking Proposition~\ref{prop. thickening P}) we obtain a homotopy commuting diagram of Kan complexes---or, equivalently, a diagram in the $\infty$-category $\Top$---as follows: 
	\eqn\label{eqn. the diagram of spaces we want to be a pullback}
	\xymatrix{
	\hom_{\MMfld^{\dd}}(P,\bullet) \ar[r] \ar[d]
		& \hom_{\MMfld^{\dd}}([-2,-1] \times X ,\bullet) \ar[d] \\
	\hom_{\MMfld^{\dd}}([1,2] \times Y,\bullet) \ar[r]
		& \hom_{\MMfld^{\dd}}(I_s \times W,\bullet).
	}
	\eqnd

\begin{prop}\label{prop. reduction to pullback of spaces}
Suppose that for every test object $\bullet$, the diagram \eqref{eqn. the diagram of spaces we want to be a pullback} is a homotopy pullback diagram of Kan complexes. Then $P$ is a pushout of~\eqref{eqn. iX iY} in the $\infty$-category $\mfld^{\dd}$.
\end{prop}

\begin{proof}
If the hypothesis is satisfied, then the image of~\eqref{eqn. first pushout} in $\mfld^{\dd}$ is a pushout diagram. (This follows, for example, from Proposition~7.4.5.13 (03BJ) of~\cite{lurie-kerodon}. In the notation of loc. cit., we take the $\infty$-category $\cC$ to be the opposite of the nerve of $\MMfld^{\dd}$.) On the other hand, we note that the compositions 
	\eqnn
	\xymatrix{
	I_s \ar[r]^-{h_1} 
		&	[1,2] \into  [-2,2] \ar[r]^-{\sigma}
		& I_s,
	}
	\qquad
	\xymatrix{
	I_s \ar[r]^-{h_{-1}} 
		&	[-2,-1] \into  [-2,2] \ar[r]^-{\sigma}
		& I_s
	}
	\eqnd
are both isotopic to $\id_{I_s}$. Thus, in $\mfld^{\dd}$ the diagram~\eqref{eqn. iX iY} is equivalent to the diagram
	\eqnn
	\xymatrix{
	[-2,-1] \times X
		&& \ar[ll]_-{h_{-1} \times i_X} I_s \times W \ar[rr]^-{h_1 \times i_Y}
		&& [1,2] \times Y
	}
	\eqnd
so $P$ is a pushout of~\eqref{eqn. iX iY}.
\end{proof}

\subsection{A convenient model for the  homotopy pullback of spaces}
\begin{notation}[$\HH^{(a)}$ and $f_X$, $f_{W,t}$, $f_Y$]
Let us suppose $Z$ is a test object with $\dim Z = \dim X$. (One may always assume this by thickening either $X$ or $Z$.) For every integer $a \geq 0$, we let
	\eqnn
	\HH^{(a)}
	\eqnd
denote the simplicial set where a $k$-simplex $f$ in $\HH^{(a)}$ is a triplet
	\eqnn
	(f_X, \{f_{W,t}\}_{t \in [-2,2]}, f_Y)
	\eqnd
of maps where
\begin{itemize}
	\item $f_X: \Delta^k \times (X \times I_s^{a}) \to \Delta^k \times (Z \times I_s^{a}) $ is a (codimension zero) smooth embedding respecting the projections to $\Delta^k$, 
	\item $f_Y: \Delta^k \times (Y \times I_s^{a}) \to \Delta^k \times (Z \times I_s^{a}) $ is a (codimension zero) smooth embedding respecting the projections to $\Delta^k$, and
	\item $\{f_{W,t}\}_{t \in [-2,2]}$ is a $t$-parametrized isotopy of embeddings, which we may think of as a smooth map
		\eqnn
		[-2,2]_t  \times \Delta^k \times (W \times I_s^{a})
		\to
		\Delta^k \times (Z \times I_s^{a})
		\eqnd
respecting the projection to $\Delta^k$, and each of whose restrictions to $t \in [-2,2]$ is a codimension-zero smooth embedding.
\end{itemize}
and the triplet must satisfy the following:
\enum
	\item $\{f_{W,t}\}$ is collared in the $t$-variable in the same sense as $h_t$ in Choice~\ref{choice. h_t}: 
		\eqn\label{eqn. collaring for fW}
		f_{W,t} = f_{W,-1} \text{ for } t \leq -1 \qquad\text{and}\qquad
		f_{W,t} = f_{W,1} \text{ for } t \geq 1.
		\eqnd
	\item For $t \leq -1$, and above every point of $\Delta^k$, the composition
		\eqnn
		\xymatrix{
		W \times I_s^{a} \ar[rr]^-{i_X \times \id} 
			&&  X  \times I_s^{a} \ar[r]^-{f_X}
			& Z \times I_s^{a} 
		}
		\eqnd
	is equal to the embedding $f_{W,t}$ (evaluated above the point of $\Delta^k$). Note that by the collaring condition, this condition need only be checked at $t=-1$.	
	\item For $t \geq 1$, and above every point of $\Delta^k$,  the composition
		\eqnn
		\xymatrix{
		W \times I_s^{a} \ar[rr]^-{i_Y \times \id} 
			&&  Y  \times I_s^{a} \ar[r]^-{f_Y}
			& Z \times I_s^{a} 
		}
		\eqnd
	is equal to the embedding $f_{W,t}$ (evaluated above the point of $\Delta^k$). Note that by the collaring condition, this condition need only be checked at $t=1$.
\enumd
\end{notation}

\begin{remark}
\label{remark. dropping simplices}
Because all constructions will respect the simplicial maps between the $\Delta^k$, we will often try to declutter notation by omitting the $\Delta^k$ variable. In practice, such notation should be interpreted to mean that---given a map $\Delta^k \times A \to \Delta^k \times B$ respecting the projections to $\Delta^k$---we are testing a condition on the induced maps $A \to B$ for every element of $\Delta^k$.
\end{remark}

\begin{notation}[$\HH^{(\infty)}$]
\label{notation. H thickening}
There are natural thickening maps $\HH^{(a)} \to \HH^{(a+1)}$ given by taking the direct product of $f_X, f_{W,t}, f_Y$ with $\id_{I_s}$. We let
	\eqnn
	\HH^{(\infty)} := \colim_{a \to \infty} \HH^{(a)}
	\eqnd
denote the increasing union. Note that because the thickening maps are cofibrations of simplicial sets, this increasing union models the homotopy colimit of simplicial sets.
\end{notation}

\begin{prop}
\label{prop. Ha is a pullback}
$\HH^{(a)}$ admits a map to the homotopy pullback of the diagram of Kan complexes
	\eqn\label{eqn. pullback of a th mapping spaces}
	\xymatrix{
		& 
	\hom_{\MMfld_{\dim X+a}}( X  \times I_s^a, Z \times  I_s^a) \ar[d] \\
	\hom_{\MMfld_{\dim X+a}}( Y \times I_s^a, Z \times   I_s^a)\ar[r]
		& \hom_{\MMfld_{\dim X+a}}( W \times I_s^a, Z \times I_s^a)
	}
	\eqnd
commuting with thickening, and this map is a homotopy equivalence. 

In particular, $\HH^{(\infty)}$ is a homotopy pullback of the diagram of Kan complexes
	\eqn\label{eqn. pullback of mfld dd spaces}
	\xymatrix{
		& 
	\hom_{\MMfld^{\dd}}( X , Z ) \ar[d] \\
	\hom_{\MMfld^{\dd}}( Y , Z )\ar[r]
		& \hom_{\MMfld^{\dd}}( W , Z )	.
	}
	\eqnd
\end{prop}

\begin{proof}
The standard model for the desired homotopy pullback is a simplicial set of triplets $(f_X,F_W,f_Y)$ where $F_W$ is a map $\Delta^1 \to \hom_{\MMfld_{\dim X+a}}( W \times I_s^a, Z \times I_s^a)$, and the evaluation of $F_W$ at $\{0\},\{1\} \in \Delta^1$ restrict to $f_X$ and $f_Y$ along $h_{-1} \times i_X$ and $h_{1} \times i_Y$, respectively.

Choosing an orientation-preserving diffeomorphism $\Delta^1 \cong [-2,2]_t$ and smoothly retracting the neighborhood $[-2,-1]$ to $\{-2\}$ and $[1,2]$ to $\{2\}$, we obtain a homotopy equivalence from the simplicial set of $\{f_{W,t}\}_{t \in [-2,2]}$ to the simplicial set of $F_W$, in a way respecting the evaluation maps at $\{0,1\} \subset \Delta^1$. Thus we have a homotopy equivalence from $\HH^{(a)}$ to the simplicial set of triplets $(f_X,F_W,f_Y)$.

The homotopy equivalence from the previous paragraph respects thickening, so we have a map of sequential diagrams from the sequence $\HH^{(0)} \to \HH^{(1)} \to \HH^{(2)} \to \ldots $ to the sequence obtained by taking the homotopy pullbacks of~\eqref{eqn. pullback of a th mapping spaces}. This map of sequences is a homotopy equivalence at the $a$th stage for every $a$, so induces a homotopy equivalence of the homotopy colimits. 

On the other hand,
\begin{itemize}
	\item The sequential colimit of~\eqref{eqn. pullback of a th mapping spaces} (as $a \to \infty$) is precisely the diagram~\eqref{eqn. pullback of mfld dd spaces}. This is in fact a homotopy colimit because the thickening maps induce cofibrations of mapping spaces.
	\item Sequential (in fact, filtered) homotopy colimits of Kan complexes commute with homotopy pullbacks (see Remark~\ref{remark. filtered colims and homotopy pullbacks}).
\end{itemize}
Thus the induced map from $\HH^{(\infty)}$ to the homotopy pullback of~\eqref{eqn. pullback of mfld dd spaces} is a homotopy equivalence of Kan complexes, proving the claim.
\end{proof}

\begin{remark}[Sequential homotopy colimits commute with homotopy pullbacks]
\label{remark. filtered colims and homotopy pullbacks}
In $\Top$, filtered homotopy colimits commute with finite homotopy limits. (This is not true of all $\infty$-categories, of course.) While this is a well-known fact, we will save the reader some trouble and explain why this is true in the special case of a sequential homotopy colimit and homotopy pullbacks. (This is the case we need in the proof of Proposition~\ref{prop. Ha is a pullback}.) 

Given a sequential diagram $(A_i) = (A_i)_{i \in I}$ (in our setting, $I$ is the linearly ordered set of natural numbers) of simplicial sets, the homotopy colimit may be computed by replacing $(A_i)$ by a projectively cofibrant diagram $(\AA_i)$, and computing the honest colimit (e.g., increasing union) of the $\AA_i$. This replacement can be made functorially, in that a map of diagrams $(A_i) \to (B_i)$ results in a map of replacements $(\AA_i) \to (\BB_i)$ while respecting compositions of maps of diagrams-- see Construction~7.5.6.8 (03CC) of~\cite{lurie-kerodon}.

(Here, by {\em replacement} we mean the data of weak homotopy equivalences $\AA_i \to A_i$ for which the two compositions $\AA_i \to A_i \to A_j$ and $\AA_i \to \AA_j \to A_j$  are equal. By a projectively cofibrant sequential diagram, we mean a diagram such that, for all $i<j$, $\AA_i \to \AA_j$ is a cofibration---i.e., monomorphism, i.e., for all $k$, injections on the set of $k$-simplices. See Example~7.5.6.4 (03C7) of~\cite{lurie-kerodon}. In particular, if $A_i$ already consists of cofibrations, the homotopy colimit is computed as the honest, point-set colimit.)

Recall that a homotopy pullback of a diagram $A \to B \leftarrow C$ of simplicial sets (in the Quillen model structure) is computed as follows: One replaces $A \to B$ by a weak equivalence followed by a Kan fibration $A \xrightarrow{\sim} A' \to B$, and one then computes the point-set pullback of the diagram $A' \to B \leftarrow C$ of Kan complexes. Up to weak homotopy equivalence, this point-set pullback is independent of the choice of the factorization $A \to A' \to B$. Even better, one can arrange for a factorization for simplicial sets for which (i) The factorization is functorial, meaning commutativity of the lefthand diagram below guarantees the commutativity of the righthand diagram
	\eqnn
	\xymatrix{
	A_i \ar[r]  \ar[d] & B_i \ar[d] \\
	A_j \ar[r] & B_j \\
	}
	\implies
	\xymatrix{
	A_i \ar[r]^{\sim}  \ar[d] & A_i' \ar[r]  \ar[d] & B_i \ar[d] \\
	A_j \ar[r]^{\sim}   & A_j'  \ar[r] &  B_j ,
	}
	\eqnd
and (ii) cofibrations are preserved, in that if $A_i \to A_j$ is a cofibration, so is the map $A_i' \to A_j'$ above. 

So given diagrams $(A_i \to B_i \leftarrow C_i)_{i \in I}$ indexed by a sequential diagram $I$, let $(\AA_i \to \BB_i \leftarrow \CC_i)_{i \in I}$ denote the sequence of diagrams obtained by functorially replacing the diagrams $(A_i),(B_i),(C_i)$ by projectively cofibrant diagrams. Because homotopy fiber products are preserved under homotopy equivalence (Corollary~8.13 of~\cite{goerss-jardine}; alternatively, see Remark 7.5.1.3 (03A2) of~\cite{lurie-kerodon}), the natural-in-$i$ maps
	\eqnn
	\holim (\AA_i \to \BB_i \leftarrow \CC_i)
	\to
	\holim (A_i \to B_i \leftarrow C_i)
	\eqnd
are all weak homotopy equivalences. Now let us functorially replace the maps $\AA_i \to \BB_i$ by fibrations $\AA_i' \to \BB_i$, so that the natural-in-$i$ maps
	\eqnn
	\holim (\AA_i \to \BB_i \leftarrow \CC_i)
	\to
	\lim(\AA_i' \to \BB_i \leftarrow \CC_i)
	\eqnd
are weak homotopy equivalences. Now contemplate the induced arrow
	\eqn\label{eqn. colim lim of simplicial sets}
	\colim_{i \in I} \lim(\AA_i' \to \BB_i \leftarrow \CC_i)
	\to
	\lim( (\colim_{i \in I} \AA_i') \to \colim_{i \in I} (\BB_i) \leftarrow \colim_{i \in I} )\CC_i)).
	\eqnd
Limits and colimits of simplicial sets are computed level-wise---for example, fiber products are given on $k$-simplices as $(A' \times_B C)_k = A'_k \times_{B_k} C_k$. It is of course a classical fact that filtered colimits commute with finite limits in sets, so we may conclude sequential colimits and fiber products commute in simplicial sets. In other words, the above arrow ~\eqref{eqn. colim lim of simplicial sets} is an isomorphism. 

On the other hand, by (ii) above, the maps $\AA_i' \to \AA_j'$ are cofibrations, which means they are levelwise injections---in particular, the maps of fiber products $\AA_i' \times_{\BB_i} \CC_i \to \AA_j' \times_{\BB_j} \CC_j$ are all level-wise injections. It follows that the colimits in~\eqref{eqn. colim lim of simplicial sets} all compute homotopy colimits. On the other hand, we also know that the arrow $\colim_{i \in I} \AA_i' \to \colim_{i \in I} \BB_i$ is a fibration because fibrations of simplicial sets are preserved under increasing unions. So all the fiber products in \eqref{eqn. colim lim of simplicial sets} are homotopy fiber products.
	
	This proves the claim that sequential  homotopy colimits and homotopy fiber products commute for simplicial sets.
\end{remark}

\begin{notation}[$r$]
For every $a \geq 0$ and every manifold $Z$ with $\dim Z = \dim P - 1$, there is a natural map
	\eqn\label{eqn. map from hom P to H}
	r: \hom_{\MMfld_{\dim P +a}}(P \times I_s^a, Z \times I_s^{a+1})
	\to
	\HH^{(a+1)}
	\eqnd
which, given an embedding $j: P \times I_s^a \to Z \times I_s^{a+1}$ (or such a collection smoothly indexed by $\Delta^k$) outputs the triplet obtained by pre-composing $j$ with the maps
	\begin{itemize}
	\item $
		\xymatrix{
		X \times I_s^{a+1} \cong (I_s \times X) \times I_s^a \ar[rr]^-{(h_{-1} \times \id) \times \id}
			&& ([-2,-1] \times X) \times I_s^a  \ar[rr]^-{\eqref{eqn. maps from X W Y into P} \times \id}
			&& P \times I_s^a
		},
		$
	\item $\{\xymatrix{
		W \times I_s^{a+1} \cong (I_s \times W) \times I_s^a  \ar[rr]^-{(h_t \times \id) \times \id}
		&& ([-2,2] \times W) \times I_s^a  \ar[rr]^-{\eqref{eqn. maps from X W Y into P} \times \id}
		&& P \times I_s^a
	}\}_{t \in [-2,2]}$, and
	\item  $
		\xymatrix{
		Y \times I_s^{a+1} \cong (I_s \times Y) \times I_s^a \ar[rr]^-{(h_{1} \times \id) \times \id}
			&& ([1,2] \times Y) \times I_s^a  \ar[rr]^-{\eqref{eqn. maps from X W Y into P} \times \id}
			&& P \times I_s^a
		}.
		$
	\end{itemize}
Informally, $r$ has the effect of substituting the $(u,x)$, $(u,w)$, $(u,y)$ variables with 
$(\alpha\beta(x), h_u(s_{a+1}))$, 
$(\alpha\beta(w), h_u(s_{a+1}))$, 
and
$(\alpha\beta(y), h_u(s_{a+1}))$, 
respectively. 
\end{notation}

\begin{remark}
\label{remark. pedantic labeling}
Let us be punctilious about the first isomorphism in each of the maps above. Write $I_s^{a+1} = I_s^{\{1,2,\ldots,a,a+1\}}$ as the set of functions from the ordered set $\{1,\ldots,a+1\}$ to $I_s$. So for example, we have the isomorphisms
	\eqnn
	\xymatrix{
	I^{a+1}_s
	=I^{\{1,\ldots,a+1\}}_s
	\cong I^{\{1,\ldots,a\}}_s \times I^{\{a+1\}}.
	}
	\eqnd
This allows us to write
	\eqnn
	X \times I^{a+1}_s 
	\cong X \times I^a_s \times I_s^{\{a+1\}} 
	\cong (X \times I_s^{\{a+1\}}) \times I^a_s
	\cong (I^{\{a+1\}}_s \times X) \times I^a_s
	\cong (I_s \times X) \times I^a_s
	\eqnd
(and like wise for $W$ and $Y$). This pedantry makes it clear that $r$ is {\em not} compatible with thickening using our convention that thickening occurs by multiplying by $I_s$ on the right.\end{remark}

\begin{remark}
\label{remark. P is pushout if map to H is equivalence}
In $\mfld^{\dd}$ there are particular equivalences 
	\eqnn
	[-2,-1] \times X \simeq X
	\qquad
	[1,2] \times Y \simeq Y,
	\qquad
	I_s \times W \simeq W
	\eqnd
given by choosing thickenings of $X,Y,W$ and permuting coordinates. Choosing such equivalences induces the map~\eqref{eqn. arrow 2} in the composition below; it is a homotopy equivalence of Kan complexes:
	\begin{align}
	\hom_{\MMfld^{\dd}}(P,Z)
	& \to
		\hom_{\MMfld^{\dd}}([-2,-1] \times X ,Z)
		\times_{\hom_{\MMfld^{\dd}}(I_s \times W,Z)}^{h}
		\hom_{\MMfld^{\dd}}([1,2] \times Y,Z) \label{eqn. arrow 1}\\
	& \to
		\holim\eqref{eqn. pullback of mfld dd spaces}\label{eqn. arrow 2} \\
	& \xrightarrow{\text{Prop~\ref{prop. Ha is a pullback}}}
		\HH^{(\infty)}. \nonumber
	\end{align}
We thus see~\eqref{eqn. arrow 1} is a homotopy equivalence if and only if the composition of the above maps in a homotopy equivalence. So consider the composition
	\eqn
	\label{eqn. map from P to H infty}
	\hom_{\MMfld^{\dd}}(P,Z) \to \HH^{(\infty)}.
	\eqnd
\end{remark}

Our goal is to now prove:

\begin{lemma}
\label{lemma. hom to H is equivalence}
For every $Z \in \ob \MMfld^{\dd}$, the map~\eqref{eqn. map from P to H infty} is a homotopy equivalence of Kan complexes.
\end{lemma}

\subsection{Proof of Lemma~\ref{lemma. hom to H is equivalence}}

\begin{choice}[$\overline{\sigma}$]
\label{choice. sigma bar}
We fix an orientation{\em-reversing} diffeomorphism
	\eqn\label{eqn. sigma bar}
	\overline{\sigma}: [-2,2] \to I_s.
	\eqnd
\end{choice}

\begin{remark}
In Choice~\ref{choice. sigma bar} we demand $\overline{\sigma}$ to be orientation-reversing. This demand is not strictly necessary; it merely allows us to simplify the proof of Proposition~\ref{prop. qr is thickening}. The trade-off is that we then need to show that $rqrq$, rather than just $rq$, is homotopic to a thickening map in Proposition~\ref{prop. rqrq is thickening}. 
\end{remark}

\begin{prop}
\label{prop. there is an isotopy from h sigmabar to id}
The map
	\eqn\label{eqn. the two variable map u and s}
	[-2,2] \times I_s \to [-2,2] \times I_s,
	\qquad
	(u,s) \mapsto (h_{u}(s),\overline{\sigma}(u))
	\eqnd
is isotopic to the identity. In fact, one can choose an isotopy of~\eqref{eqn. the two variable map u and s} to the identity through maps for which the locus $\{u \leq -1\}$ has image contained in itself, and likewise for the locus $\{u \geq 1\}$. 
\end{prop}

\begin{figure}[ht]
    \begin{equation}\nonumber
			\xy
			\xyimport(8,8)(0,0){\includegraphics[width=6in]{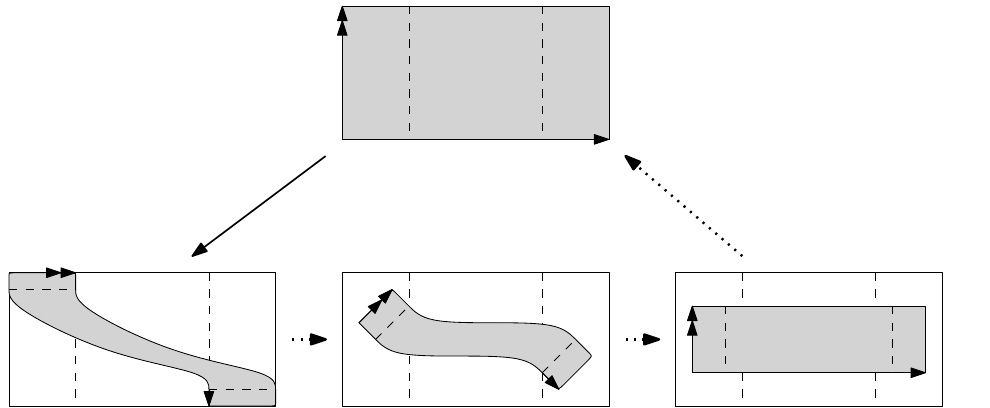}}
				,(1.9,4.4)*{\eqref{eqn. the two variable map u and s}};
			\endxy
    \end{equation}
    \caption{
    An isotopy of the map~\eqref{eqn. the two variable map u and s} to the identity. The horizontal direction is the $[-2,2]_u$ direction, while the vertical direction is the $I_s$ direction, with orientations indicated (in the directions of increasing values of $u$ and of $s$). The dotted arrows indicate isotopies, with the images of each stage of the isotopy indicated as the shaded regions. Note also that the dashed lines indicate the lines $u = \pm 1$ and their images under various stages of the isotopy.
    }
    \label{figure.u-s-isotopy}
\end{figure}

\begin{proof}
This follows from our assumptions on $h$ (Choice~\ref{choice. h_t}). See Figure~\ref{figure.u-s-isotopy}.
\end{proof}

\begin{notation}[$q'$]
\label{notation. q'}
Fix $a \geq 1$, and fix a triplet $(f = f_X, \{f_{W,t}\}_{t \in [-2,2]}, f_Y) \in \HH^{(a)}$. We define an embedding
	\eqnn
	q'(f): P' \times I_s^a \to (Z \times I_s^{a}) \times I_s
	\eqnd
to be the union of the three compositions
	\begin{align}
		\xymatrix{
		[-2,-1] \times (X  \times I_s^{a})
		\cong
		(X  \times I_s^{a}) \times [-2,-1]
		\ar[r]^-{\id \times \overline{\sigma}}
		& (X \times I_s^{a}) \times I_s
		\ar[rr]^{f_X \times \id}
		&& (Z \times I_s^{a}) \times I_s.
		}\label{eqn. r on X}
		\\
		\xymatrix{
		[-1,1] \times (W \times I_s^a)
		\cong (W \times I_s^a) \times [-1,1]
		\ar[r]^-{\id \times \overline{\sigma}}
			& (W \times I_s^a) \times I_s
		\ar[rr]^{(f_{W,\overline{\sigma}^{-1}},\id)}
		&&  (Z \times I_s^a) \times I_s
		}\label{eqn. r on W} \\
		\xymatrix{
		[1,2] \times (Y  \times I_s^{a})
		\cong
		(Y  \times I_s^{a}) \times [-2,-1]
		\ar[r]^-{\id \times \overline{\sigma}}
		& (Y \times I_s^{a}) \times I_s
		\ar[rr]^{f_Y \times \id}
		&& (Z \times I_s^{a}) \times I_s.
		}\label{eqn. r on Y}
	\end{align}
Note that we may define a map out of $P'$ as a union of these three maps by definition of $P'$~\eqref{eqn. P before smoothing}. That the union of these maps is smooth follows from the collaring condition~\eqref{eqn. collaring for fW}. That the union is a diffeomorphism follows from the facts that $\overline{\sigma}$ is a diffeomorphism~\eqref{eqn. sigma bar} and that the maps $f_X, f_Y$ and $f_{W,t}$ for each $t$, are embeddings by definition of $\HH^{(a)}$.
\end{notation}

\begin{example}
On a triplet $(u,w, \vec s) \in [-1,1] \times (W \times I^a_s)$, \eqref{eqn. r on W} evaluates as
	\eqnn
	(f_{W,u}(w,\vec s),\overline{\sigma}(u)) \in (Z \times I^a_s) \times I_s.
	\eqnd
\end{example}

\begin{notation}[$q$]
\label{notation. q}
We have defined $q'(f)$ as a map with codomain $P' \times I_s^a$. By pre-composing $q'(f)$ with the map $\alpha: P \to P'$~\eqref{eqn. alpha from P to P'}, we obtain a map of simplicial sets
	\eqnn
	q: 
	\HH^{(a)}
	\to \hom_{\MMfld_{\dim P + a}}(P \times I_s^a, (Z \times I_s^a) \times I_s),
	\qquad
	f \mapsto (q'(f)) \circ \alpha .
	\eqnd
(To define $q$ and $q'$ on a $k$-simplex $f$, we follow Remark~\ref{remark. dropping simplices} to note that \eqref{eqn. r on X},~\eqref{eqn. r on W}, and~\eqref{eqn. r on Y} all make sense as parametrizing functions parametrized over points of the $k$-simplex.)
\end{notation}

We now have a sequence of Kan complexes
	\eqn\label{eqn. q r sequence}
	\xymatrix{
	\vdots \ar[d]^r \\
	\HH^{(a)} \ar[d]^{q} \\
	\hom_{\MMfld_{\dim P +a}}(P \times I_s^a, Z \times I_s^{a+1}) \ar[d]^{r} \\
	\HH^{(a+1)} \ar[d]^{q} \\
	\hom_{\MMfld_{\dim P +a+1}}(P \times I_s^{a+1}, Z \times I_s^{a+2}) \ar[d]^{r} \\
	\vdots.
	}
	\eqnd

\begin{prop}
\label{prop. qr is thickening}
For every $a \geq 0$, the composition $qr$ is homotopic to the thickening map. (That is, $qr$ is homotopic to the map of Kan complexes $\hom_{\MMfld_{\dim P +a}}(P \times I_s^a, Z \times I_s^{a+1}) \to \hom_{\MMfld_{\dim P +a+1}}(P \times I_s^{a+1}, Z \times I_s^{a+2})$---defined as in~\eqref{eqn. thickening functor} by taking the direct products of domain and codomain manifolds with $I_s = [0,1]$---see Notation~\ref{notation. I_s}.)
\end{prop}

\begin{proof}[Proof of Proposition~\ref{prop. qr is thickening}.]
Fix an embedding $j': P' \times I_s^{a-1} \to Z \times I_s^a$. 
Consider the restriction of $q'r(j'\circ \alpha )$ along the map 
	\eqn\label{eqn. map from W thickened into P}
	([-2,2] \times W \times I^{a-1}_s) \times I_s
	\xrightarrow{(\eqref{eqn. the maps we want to be smooth} \times \id) \times \id}
	(P' \times I^{a-1}_s) \times I_s.
	\eqnd
Parsing through the definitions---\eqref{eqn. map from hom P to H},~\eqref{eqn. r on X},~\eqref{eqn. r on W}, and~\eqref{eqn. r on Y}---we see this restriction is obtained by the composition
	\begin{align}
		[-2,2] \times (W \times I^{a}_s)
		 	& \cong (W  \times I_s^{a}) \times [-2,2] \label{eqn. the isotopy we want to make trivial 1}  \\ 
			& \xrightarrow{\id \times \overline{\sigma}} (W \times I_s^{a}) \times I_t \nonumber\\
		& \cong  \left( (I_s^{\{a\}} \times W) \times I_s^{\{1,\ldots,a-1\}}\right) \times I_t \label{eqn. the isotopy we want to make trivial pedantic label} \\
		& \xrightarrow{( ((h_{\overline{\sigma}^{-1}(t)} \times \id) \times \id) , \id)} 
			 \left(([-2,2] \times W) \times I^{a-1}_s\right) \times I_t \label{eqn. the isotopy we want to make trivial 2} \\
		& \xrightarrow{(\eqref{eqn. the maps we want to be smooth} \times \id) \times \id} 
			 \left(P' \times I^{a-1}_s\right) \times I_t \label{eqn. glue together to id of P'} \\
		& \xrightarrow{(\beta \times \id) \times \id} 
			 \left(P \times I^{a-1}_s\right) \times I_t \nonumber \\
		& \xrightarrow{(j'\circ (\alpha \times \id)) \times \id} 
			 (Z \times I^{a-1}_s) \times I_t. \nonumber 
	\end{align}
In the above notation, the intervals $I_s = I_t$ are the same set. The subscripts indicate that (in the relevant lines) we use $t$ to denote a variable representing the elements of the interval $I_t$, to remove the ambiguity regarding the copy of $I$ on which $h_{\overline{\sigma}^{-1}(t)}$ depends. Having disambiguated, we will henceforth revert to writing $I_s$ for all thickening intervals. 

Note also that in line~\eqref{eqn. the isotopy we want to make trivial pedantic label}, we have used notation from Remark~\ref{remark. pedantic labeling}.

In fact, by replacing all instances of $W$ by $X$, and all instances of $[-2,2]$ by $[-2,-1]$---which has the effect of turning $h_{\sigma^{-1}(t)}$ into $h_{-1}$---the above lines also describe the restriction of $q'r(j'\circ \alpha )$ along 
	\eqnn
	([-2,-1] \times X \times I^{a-1}_s ) \times I_s
	\xrightarrow{(\eqref{eqn. the maps we want to be smooth} \times \id) \times \id}
	(P' \times I^{a-1}_s) \times I_s,
	\eqnd
and likewise for $Y$ (by replacing the intervals with $[1,2]$). 

Now consider the compositions from line~\eqref{eqn. the isotopy we want to make trivial 1} through line~\eqref{eqn. the isotopy we want to make trivial 2} for each of $X,W,Y$ (and their relevant intervals $[-2,-1]$, $[-1,1]$ and $[1,2]$, respectively). The compositions define self-embeddings 
	\eqnn
	\eta_X,
	\qquad
	\eta_W,
	\qquad
	\eta_Y,
	\eqnd
of the manifolds
	\eqnn
	([-2,-1] \times X \times I^{a-1}_s ) \times I_s,
	\qquad
	([-2,2] \times W \times I^{a-1}_s ) \times I_s,
	\qquad
	([1,2] \times Y \times I^{a-1}_s ) \times I_s.
	\eqnd
Each of these self-embeddings is a direct product with the identity embeddings of $X,W,Y$, respectively. In particular, the maps $\eta_X,\eta_W,\eta_Y$ glue together to form a single embedding
	\eqnn
	\eta:  \left(P' \times I^{a-1}_s\right) \times I_s \to  \left(P' \times I^{a-1}_s\right) \times I_s,
	\eqnd
On the interval components---i.e., ignoring the $P'$ factor---$\eta$ acts as
	\eqnn
	[-2,2] \times I^a_s
		\xrightarrow{\cong} 
		[-2,2] \times I_s^{a-1} \times I,
		\qquad
	(u, s_1,\ldots,s_a) 
		\mapsto 
		(h_{u}(s_a),( s_1,\ldots,s_{a-1}) , \overline{\sigma}(u))). 
	\eqnd
In other words, the composition acts trivially on all but two factors: It acts on the $u$ and $s_a$ components precisely by the map~\eqref{eqn. the two variable map u and s} from Proposition~\ref{prop. there is an isotopy from h sigmabar to id}. So the isotopy guaranteed by Proposition~\ref{prop. there is an isotopy from h sigmabar to id} induces isotopies
	\eqnn
	\eta_X \sim \id,
	\qquad
	\eta_W \sim \id,
	\qquad
	\eta_Y \sim \id.
	\eqnd
(Note the condition that the isotopy preserves the loci $u \leq -1, u \geq 1$ guarantees that $\eta_X$ is isotoped through $([-2,-1] \times X \times I^{a-1}_s ) \times I_s$, and likewise for $Y$.)
These isotopies glue together to exhibit an isotopy 
	\eqn\label{eqn. eta id isotopy}
	\eta \sim \id.
	\eqnd
At this point we can appreciate that gluing together the $W,X,Y$ versions of~\eqref{eqn. glue together to id of P'} yields the identity map of $P'$. This means 
	\eqnn
	q'r(j'\circ \alpha )
	=
	((j'\circ (\alpha \times \id_{I^{a-1}_s}) \circ (\beta \times \id_{I^{a-1}_s})) \times \id) \circ \eta
	\eqnd
and applying the isotopy~\eqref{eqn. eta id isotopy} yields an isotopy
	\eqnn 
	q'r(j'\circ \alpha )
	\sim
	(j'\circ (\alpha \times \id_{I^{a-1}_s}) \circ (\beta \times \id_{I^{a-1}_s})) \times \id 
	\eqnd
which, because $\alpha\beta $ is isotopic to $\id$, is isotopic to $j' \times \id$. By varying $j'$, we find that the diagram of Kan complexes
	\eqnn
	\xymatrix{
		\hom_{\MMfld_{\dim P + a - 1}}(P \times I^{a-1}_s, Z \times I_s^a)
		\ar[r]^-{r}
		& \HH^{(a)}
		\ar[r]^-{q} \ar[dr]^{q'}
		&	\hom_{\MMfld_{\dim P + a}}(P \times I^{a}_s, Z \times I_s^{a+1})
		\\
	\sing\emb'(P',\times I^{a-1}_s, Z \times I_s^a)
		\ar[u]^-{\alpha^*} \ar[rr]^-{-\times \id_{I_s}}
		&&
		\sing\emb'(P' \times I^{a}_s, Z \times I_s^{a+1}) \ar[u]^{\alpha^*}
	}
	\eqnd
commutes up to homotopy. (In fact, the triangle on the right commutes on the nose.) By noting that $\alpha^*$ is a homotopy equivalence~\eqref{eqn. emb P is emb'P'} that respects thickening up to homotopy, we are finished.
\end{proof}

\begin{notation}[$H$]
\label{notation. H self map}
We will now compute the iterate $rqrq(f)$. For this, we introduce the following function for sake of brevity:
	\eqnn
	H: I_s \times I_s \times [-2,2]_t \to [-2,2]_t,
	\qquad
	H(s_{a+1},s_{a+2},t) = h_{h_t(s_{a+2})}(s_{a+1}).
	\eqnd
\end{notation}

\begin{remark}
Note that---by the collaring condition on $h_t$ (Choice~\ref{choice. h_t})---we have
	\eqnn
	H=
	\begin{cases}
	h_{-1}(s_{a+1}) & t \in [-2,-1] \\
	h_{1}(s_{a+1}) & t \in [1,2] .
	\end{cases}
	\eqnd
\end{remark}

\begin{prop}
\label{prop. H homotopies}
The self-map of $[-2,2]_t \times I_s \times I_s$ given by
	\eqn\label{eqn. rqrq map's relevant bit}
	(t,s_{a+1},s_{a+2})
	\mapsto
	(H, \overline{\sigma}H, \overline{\sigma}h_t(s_{a+2}))
	\eqnd
is homotopic to the identity. Moreover, this homotopy $\{k_\tau\}_{\tau \in [0,1]}$ can be chosen through smooth maps 
	\eqnn
	k_\tau: [-2,2]_t \times I_s \times I_s \to [-2,2]_t \times I_s \times I_s
	\eqnd
such that
\enum[(a)]
	\item\label{item. embedding requirement of homotopy} For any value of $t$ and $\tau$, the composition
	\eqn\label{eqn. t restricted k}
	\{t\} \times I_s \times I_s \xrightarrow{k_\tau} [-2,2]_t \times I_s \times I_s \to I_s \times I_s
	\eqnd
	is a (codimension zero) smooth embedding. Here, the first arrow is the restriction of $k_\tau$ to the indicated domain, and the last map is the projection forgetting the $[-2,2]$ factor.
	\item\label{item. collaring requirement of homotopy} For every value of $\tau$, the composition~\eqref{eqn. t restricted k} is collared in the $t$-variable. This means   that for all $t \leq -1$, the projections of $k_\tau(t,-,-)$ and  $k_\tau(-1,-,-)$ to $I_s \times I_s$ are equal. Likewise, for all $t \geq 1$, the projections of $k_\tau(t,-,-)$ and $k_\tau(1,-,-)$ to $I_s \times I_s$ are equal. (However, their projections to the $t$ variable may differ.) 
	\item\label{item. X-Y preserving requirement of homotopy} For every value of $\tau$, $k_\tau$ restricts to a self-map of $[-2,-1]_t \times I_s \times I_s$---meaning the image of $[-2,-1]_t \times I_s \times I_s$ under $k_\tau$ is contained in $[-2,-1]_t \times I_s \times I_s$---and to a self-map of $[1,2]_t \times I_s \times I_s$.
\enumd
\end{prop}

\begin{proof}
Let us first focus on the first factor of~\eqref{eqn. rqrq map's relevant bit}, i.e., the function $H$ taking value in the $[-2,2]_t$ component of the image. We note that the conditions~\eqref{item. embedding requirement of homotopy} and~\eqref{item. collaring requirement of homotopy} are preserved by homotopies $\{k'_\tau\}_{\tau}$ that only affect this first factor. On the other hand, it is straightforward to produce a homotopy $\{k'_\tau\}_{\tau}$ without changing the latter two factors of~\eqref{eqn. rqrq map's relevant bit}, smoothly homotoping the first factor from $H$ to $t$, while preserving~\eqref{item. X-Y preserving requirement of homotopy}. Choose such a $\{k'_\tau\}_{\tau}$; we have homotoped~\eqref{eqn. rqrq map's relevant bit} to the map
	\eqnn
	(t,s_{a+1},s_{a+2})
	\mapsto
	(t, \overline{\sigma}H, \overline{\sigma}h_t(s_{a+2})).
	\eqnd
Now, for every fixed $t$, note it is possible to {\em isotope} the map $(\overline{\sigma}H, \overline{\sigma}h_t(s_{a+2}))$ to the identity map of $I_s \times I_s$. It is clear one can do this smoothly in $t$, so choose such isotopies smoothly in $t$, and constantly along the intervals $-2 \leq t \leq -1$ and $1 \leq t \leq 2$. Such a choice yields a smooth homotopy $\{k''_\tau\}_{\tau}$ (whose projection to the interval $[-2,2]_t$ is $\tau$-independent) preserving all the conditions required by the Proposition. 

Concatenating $\{k''_\tau\}_{\tau}$ after $\{k'_\tau\}_{\tau}$, we have obtained the desired homotopy $\{k_{\tau}\}$. 
\end{proof}

\begin{prop}
\label{prop. rqrq is thickening}
For every $a \geq 0$, the composition $rqrq$ is homotopic to the square of the thickening map---i.e., homotopic to the map $\HH^{(a)} \to \HH^{(a+2)}$ obtained by performing twice the thickening map in Notation~\ref{notation. H thickening}. 
\end{prop}

\begin{proof}[Proof of Proposition~\ref{prop. rqrq is thickening}.]
Let us fix $f \in \HH^{(a)}$ and compute the triplet $rq(f)$.  For brevity, we let
	\eqnn
	\vec{s} := (s_1,\ldots,s_a) \in I_s^a.
	\eqnd
The map $q(f)$ sends a point $(p,\vec{s}) \in P \times I_s^a$ to
	\eqnn
	\begin{cases}
	(f_{X}(x, \vec{s}), \overline{\sigma}(u))
		& \text{ if } p = (u,x )  \in [-2,-1] \times X  \\
	(f_{W,u}(w, \vec{s}), \overline{\sigma}(u))
		& \text{ if } p = (u,w ) \in [-2,2] \times W   \\
	(f_{Y}(y, \vec{s}), \overline{\sigma}(u))
		& \text{ if } p = (u,y )  \in [1,2] \times Y .
	\end{cases}
	\eqnd
(Note that for points $p$ in $[-2,-1] \times W$, the above is well-defined thanks to the collaring conditions on $h$ and on $f_{W,t}$; likewise for points in $[1,2] \times W$.) 
For every $t \in [-2,2]$, we have that
	\eqnn
	(rq(f))_{W,t} : W \times I_s^{a+1} \to (Z \times I_s^a) \times I_s
	\eqnd
is given by the  formula
	\eqnn
	(w,\vec{s},s_{a+1})
	\mapsto
	\left(
	f_{W,h_t(s_{a+1})}(\alpha_{h_t(s_{a+1})}\beta_{h_t(s_{a+1})}(w),\vec{s}),\overline{\sigma}h_t(s_{a+1})
	\right).
	\eqnd
Let us explain this notation. Clearly $\alpha$ respects the $u$-coordinate of $P'$, as does $\beta$ (Remark~\ref{remark. beta respects u coordinate}). Thus it makes sense to fix a $u$-coordinate such as $h_t(s_{a+1})$ and restrict $\alpha$ and $\beta$ as functions with domain and codomain given by fibers above this $u$-value. In fact, $\alpha_u$ is always the identity, so we may omit $\alpha$ in much of what follows.

Likewise, the maps
	\eqnn
	(rq(f))_X: X \times I_s^{a+1} \to (Z \times I_s^a) \times I_s
	\qquad\text{and}\qquad
	(rq(f))_Y: Y \times I_s^{a+1} \to (Z \times I_s^a) \times I_s
	\eqnd
have formulas
	\eqnn
	(x,\vec{s},s_{a+1})
	\mapsto
	\left(
	f_{X}(\beta_{h_{-1}(s_{a+1})}(x),\vec{s}),\overline{\sigma}h_{-1}(s_{a+1})
	\right)
	\eqnd
and
	\eqnn
	(y,\vec{s},s_{a+1})
	\mapsto
	\left(
	f_{Y}(\beta_{h_{1}(s_{a+1})}(y),\vec{s}),\overline{\sigma}h_{1}(s_{a+1})
	\right),
	\eqnd
respectively. Note that in writing the formula for $(rq(f))_X$, we have used that $h_{-1}(s_{a+1}) \in [-2,-1]$ and that $f_{W,t} = f_{W,-1} = f_X \circ i_W$ for $t \in [-2,-1]$; and likewise for the formula for $(rq(f))_Y$.

Now that we have computed $rq$, we may iterate to compute the following.
\begin{itemize}
	\item $rqrq(f)_{X}$ sends an element $(x,\vec{s},s_{a+1},s_{a+2}) \in X \times I^{a+2}_s$ to
    	\eqnn
    	\left(
    	f_X( \beta_{h_{-1}(s_{a+1})} \beta_{h_{-1}(s_{a+2})}(x),\vec s), \overline{\sigma}h_{-1}(s_{a+1}), \overline{\sigma}h_{-1}(s_{a+2}) 
    	\right),
    	\eqnd
	\item $rqrq(f)_{Y}$ sends an element $(y,\vec{s},s_{a+1},s_{a+2}) \in Y \times I^{a+2}_s$ to
    	\eqnn
    	\left(
    	f_Y(\beta_{h_1(s_{a+1})}\beta_{h_1(s_{a+2})}(y),\vec{s}),\overline{\sigma}h_1(s_{a+1}),\overline{\sigma}h_1(s_{a+2})
    	\right),
    	\eqnd
	\item and for all $t \in [-2,2]$, $rqrq(f)_{W,t}$ sends an element $(w,\vec{s},s_{a+1},s_{a+2})$ to
    	\eqnn
    	\left(
    	f_{W,H}(\beta_H \beta_{h_t(s_{a+2})}(w),\vec{s}),\overline{\sigma}H,\overline{\sigma}h_t(s_{a+2})
    	\right)
    	\in
    	(Z \times I_s^a) \times I_s \times I_s.
    	\eqnd
\end{itemize}
Thus, each of these maps can be understood as a composition
	\eqn\label{eqn. rqrq as a composition}
	\xymatrix{
	[-2,2]_t \times \bullet \times I^a_s \times I_s \times I_s
	\ar[r]^{\eqref{eqn. a map with too many betas}}
		& ([-2,2]_t \times \bullet \times I^a_s) \times I_s \times I_s
	\ar[rr]^-{(f_{\bullet,t}) \times \id \times \id}
		&& (Z \times I_s^a) \times I_s \times I_s
	}
	\eqnd
where the first map is given in coordinates by
	\eqn\label{eqn. a map with too many betas}
	(t,\bullet,\vec{s},s_{a+1},s_{a+1})
	\mapsto 
	\left(
		H, \beta_H \beta_{h_t(s_{a+1})}(\bullet), \vec s, \overline{\sigma}H,\overline{\sigma}h_t(s_{a+2})
	\right).
	\eqnd
Here, the $\bullet$ is a stand-in for $X,W,Y$ (or $x,w,y$ in the coordinate formula), it is understood that $f_{X,t} = f_X$ and $f_{Y,t} = f_Y$, and we only evaluate the first map on those elements of $[-2,2] \times \bullet$ that are elements of $P'$. 

Because $\beta_u$ is a flow by a vector field (hence isotopic to the identity) we may find a path in $\HH^{(a+2)}$ from~\eqref{eqn. a map with too many betas} to the map
	\eqnn
	(t,\bullet,\vec{s},s_{a+1},s_{a+1})
	\mapsto 
	\left(
		H, \bullet, \vec s, \overline{\sigma}H,\overline{\sigma}h_t(s_{a+2})
	\right).
	\eqnd
By Proposition~\ref{prop. H homotopies}, this map enjoys a path in $\HH^{(a+2)}$ to the (identity) map
	\eqnn
	(t,\bullet,\vec{s},s_{a+1},s_{a+1})
	\mapsto 
	\left(
		t, \bullet, \vec s, s_{a+1}, s_{a+2}
	\right).
	\eqnd
Thus, post-composing the resulting isotopy from \eqref{eqn. a map with too many betas} to $\id$ with $(f_{\bullet,t}) \times \id \times \id$, we witness an isotopy from~\eqref{eqn. rqrq as a composition} to a two-fold thickening of $f \in \HH^{(a)}$. Because the isotopy between $rqrqf$ and $f \times \id \times \id$ is witnessed by pre-composing $f \times \id \times \id$ by a series of isotopies independent of $f$, the claim is proven.
\end{proof}

\begin{proof}[Proof of Lemma~\ref{lemma. hom to H is equivalence}.]
The sequential diagram ~\eqref{eqn. q r sequence} has a cofinal subdiagram as follows:
	\eqn\label{eqn. qr colimit diagram}
	\ldots
	\xrightarrow{qr}
	\hom_{\MMfld_{\dim P +a}}(P \times I_s^a, Z \times I_s^{a+1})
	\xrightarrow{qr}
	\hom_{\MMfld_{\dim P +a+1}}(P \times I_s^a, Z \times I_s^{a+2})
	\xrightarrow{qr}
	\ldots.
	\eqnd
By Proposition~\ref{prop. qr is thickening}, this sequence is equivalent to the sequence of thickenings, hence has colimit given by $\hom_{\mfld^{\dd}}(P,Z)$. 

There is another cofinal subdiagram of the following form:
	\eqn\label{eqn. rqrq colimit diagram}
	\ldots
	\xrightarrow{rqrq}
	\HH^{(a)}\xrightarrow{rqrq}\HH^{(a+2)}
	\xrightarrow{rqrq}
	\ldots.
	\eqnd
(In fact, there are two such---one could take $a$ to be odd or even. The choice is immaterial to us.) By Proposition~\ref{prop. rqrq is thickening}, this sequence is equivalent to the sequence with the same objects, with maps given by (double) thickenings. Thus this subdiagram has colimit given by $\HH^{(\infty)}$. 

It is now straightforward to check that the induced map on homotopy groups
	\eqnn
	\xymatrix{
		\pi_* \hom_{\mfld^{\dd}}(P,Z) \cong \pi_*\hocolim\eqref{eqn. qr colimit diagram} \ar[r]^-{\cong}
			& \pi_*\hocolim\eqref{eqn. q r sequence} \ar@{-->}[r]^-{\cong}
			& \pi_* \hocolim\eqref{eqn. rqrq colimit diagram} \cong \pi_* \HH^{(\infty)}
	}
	\eqnd
(where the dashed arrow is the inverse to the isomorphism induced by the cofinality of~\eqref{eqn. rqrq colimit diagram}) is the same map as that induced on homotopy groups by~\eqref{eqn. map from P to H infty}. 
\end{proof}

\section{Proofs of the main results}
\subsection{Finite colimits: Existence and preservation}
\label{section. finite colimits in mfld}

\begin{proof}[Proof of Theorem~\ref{theorem. pushouts exist}.]
It is clear that $\mfld^{\dd}$ has an initial object---the empty manifold---so to see $\mfld^{\dd}$ has all finite colimits, it suffices to prove that $\mfld^{\dd}$ has pushouts (Corollary 4.4.2.4 of~\cite{lurie-htt}). 

Given two morphisms as in~\eqref{eqn. iX iY}, Proposition~\ref{prop. reduction to pullback of spaces} shows that $P$ is a pushout if and only if the map~\eqref{eqn. arrow 1} is a homotopy equivalence. By Remark~\ref{remark. P is pushout if map to H is equivalence}, we are reduced to showing that~\eqref{eqn. map from P to H infty} is a homotopy equivalence.  This is Lemma~\ref{lemma. hom to H is equivalence}.
\end{proof}

\begin{remark}
\label{remark. slice categories create colimits}
Let $\cC$ be an $\infty$-category and fix an object $B \in \ob \cC$. Fix a functor $f: D \to \cC_{/B}$ to the slice $\infty$-category.
Then a diagram $D^{\triangleright} \to \cC_{/B}$ is a colimit diagram if and only if the composition $D^{\triangleright} \to \cC_{/B} \to \cC$ is a colimit diagram (see, for example, Corollary 7.1.3.20 (02KC) of~\cite{lurie-kerodon}). 
\end{remark}

	\begin{prop}\label{prop. main functor preserves finite colims}
	The functor~\eqref{eqn. main functor} preserves finite colimits.
	\end{prop}

\begin{proof}
The isotopy-commuting diagram of smooth manifolds~\eqref{eqn. first pushout}---because it consists of codimension zero embeddings---induces a homotopy-commuting diagram of frame bundles. By quotienting by the free $O$ action, we thus obtain a homotopy-commuting diagram of spaces over $BO$. (It may help, or confuse, the reader that this resulting diagram can be notated identically to~\eqref{eqn. first pushout}.) It is classical that~\eqref{eqn. first pushout} is in fact a colimit diagram in the $\infty$-category of topological spaces. (After all, $P$ models a  homotopy pushout of~\eqref{eqn. iX iY}, as $P$ is homotopy equivalent to~\eqref{eqn. P before smoothing}.) On the other hand, a diagram in the slice over-category $\TTop_{/BO}$ is a colimit diagram if and only if its image in $\TTop$ is (Remark~\ref{remark. slice categories create colimits}).

This shows that~\eqref{eqn. main functor} preserves pushouts. The functor also preserves initial objects, as it sends the empty manifold to the empty space. Thus the functor preserves all finite colimits (Corollary 4.4.2.5 of~\cite{lurie-htt}).
\end{proof}
\subsection{The point generates thickened manifolds}

\begin{prop}[Handle attachments are pushouts]\label{prop. handle attachments are pushouts}
Let $X$ be a manifold of dimension $d$, possibly with boundary, and let $X'$ be obtained from $X$ by attaching an index $k$ handle along $\del X$. Then the induced diagram $\Delta^1 \times \Delta^1 \to \mfld^{\dd}$, which we informally draw as
	\eqnn
	\xymatrix{
	S^{k-1} \times D^{d-k} \ar[r] \ar[d] & X \ar[d] \\
	D^{k} \times D^{d-k} \ar[r] & X'
	}
	\eqnd
is a pushout square. (Note that the top horizontal arrow is modeled as the codimension zero embedding of $(S^{k-1} \times D^{d-k}) \times [0,1]$ into a collar neighborhood of $\del X \subset X$.)
\end{prop}

\begin{proof}
We may assume there exists a Morse function $f: X' \to \RR$ realizing the handle attachment---concretely, $f^{-1}[b-\delta,b+\delta]$ is the attaching handle $D^k \times D^{n-k}$ inside $X'$, we assume $b$ is the lone critical value in $[b-\delta,\infty]$, and that $X = f^{-1}(-\infty,b-\delta]$. Using the induced map
	\eqnn
	\xymatrix{
	X' \times [-2,2]_v
	\ar[r]^{f \times p_v}
	& \RR \times [-2,2]_v
	}
	\eqnd
one may describe the space $P'$ from~\eqref{eqn. diagram to P} as the subspace of $X' \times [-2,2]_v$ given by the union of the three subsets
	\begin{itemize}
	\item $\{ (f(x) \leq b-{\frac \delta 2}) \, \& \, (v \leq 1\}$,
	\item $\{(f(x) \in [b-\delta,b+\delta]) \, \& \, (v \in [1-\epsilon,1+\epsilon])\}$, and
	\item $\{ (f(x) \geq b+{\frac \delta 2}) \, \& \, (v \geq 1)\}$.
	\end{itemize}
It is then clear that the inclusion of the space $P$ (Construction~\ref{construction. P}) into $X' \times [-2,2]$ is an isotopy equivalence. Because $X' \simeq X' \times [-2,2]$ in $\mfld^{\dd}$ and $P$ was already shown to be a pushout, the rest is straightforward.
\end{proof}

The following is a categorical manifestation of the geometric tautology that manifolds are made out of disks: 

	\begin{lemma}\label{lemma. point generates} 
	The point generates $\mfld^{\dd}$ under finite colimits. 
	\end{lemma}

\begin{proof}[Proof of Lemma~\ref{lemma. point generates}.]
Let $X$ be a compact manifold, possibly with corners and boundary. By Remark~\ref{remark. isotopy equivalence removes corners} we may assume $X$ has no corners, but possibly non-empty boundary. Then $X$ admits a Morse function $f: X \to \RR$ for which $X = f^{-1}[-\infty,b]$ and $f^{-1}(b) = \del X$, with $b$ a regular value. (See for example Theorem~2.5 of~\cite{milnor-1965-h-cobordism-lectures}. Using the notation there, one considers $V_0 = \emptyset$.) In the usual way, $f$ creates a filtration of $X$
	\eqnn
	X = X_N \supset X_{N-1} \supset \ldots \supset X_1 \supset X_0 \supset X_{-1} = \emptyset
	\eqnd
where each $X_i$ is obtained by attaching a single handle of index $a_i$ to $X_{i-1}$. By Proposition~\ref{prop. handle attachments are pushouts}, this realizes $X_i$ as a pushout of $X_{i-1}$, a (thickened) sphere, and a (thickened) point. On the other hand, by induction on $k$, Proposition~\ref{prop. handle attachments are pushouts} also shows that any sphere is generated in $\mfld^{\dd}$ under finite colimits by a point. 
\end{proof}

\subsection{Fully faithful}
\label{section. fully faithful}

\begin{prop}\label{prop. hom from point is frame bundle}
For any smooth, compact manifold $X$ (possibly with corners) there exists a natural homotopy equivalence of Kan complexes
	\eqn\label{eqn. emb to Fr}
	\hom_{\MMfld^{\dd}}(pt,X) \to \sing \Fr(X).	
	\eqnd
Here, $\sing\Fr(X)$ is the Kan complex of singular chains associated to the stable frame bundle $\Fr(X)$, and $\Fr$ is modeled as in Remark~\ref{remark. model for Fr}.
\end{prop}

\begin{proof}
Set $d = \dim X$. The space of smooth embeddings $[0,1]^d \to X$ fits into a fiber sequence 
			\eqnn
			GL(\RR^d, T_{\ev_0} X) \to \emb([0,1]^d, X) \xrightarrow{ev_0} X
			\eqnd
where the last map is the evaluation at the origin $0 \in [0,1]^d$, and the fiber is identified with the space of linear isomorphisms from the tangent space $\RR^d \cong T_0 [0,1]^d$ to $T_{ev_0} X$. (Note that if $X$ has boundary or corners, we may replace $X$ by its interior to obtain a homotopy-equivalent fiber sequence.) It follows that the map 
		 \eqnn
		 \emb([0,1]^d,X) \to \Fr_d(X)
		 \eqnd
to the frame bundle---given by sending an embedding $j: [0,1]^d \to X$ to the derivative-induced framing $\RR^d \cong T_0[0,1]^d \cong T_{j(0)} X$---is a map of fibrations over $X$ (or over the interior of $X$) with homotopy equivalent fibers, hence a homotopy equivalence. By taking the colimit of the induced maps of Kan complexes 
	\eqnn
	\sing\emb([0,1]^{d+k}, X \times [-1,1]^k) \to \sing \Fr_{d+k} (X)
	\eqnd
as $k \to \infty$, and noting that the transition maps from $k$ to $k+1$ are all cofibrations, we find that the induced map~\eqref{eqn. emb to Fr} is indeed a homotopy equivalence of Kan complexes.

Now, given any smooth codimension zero embedding $X \times [0,1]^k \to X' \times [0,1]^{k'}$ (note this forces $d+k = d'+k'$), the diagram of Kan complexes below commutes:
		 \eqnn
		 \xymatrix{ 
		 	&\hom_{\MMfld^{\dd}} ([0,1]^{d+k}  ,X\times [0,1]^{k} ) \ar[r] \ar[d]& \sing\Fr_{d+k}(X \times [0,1]^{k}) \ar[d] \\ 
		 	&\hom_{\MMfld^{\dd}} ([0,1]^{d'+k'}  ,X'\times [0,1]^{k'})  \ar[r]   &\sing\Fr_{d' +k'}(X' \times [0,1]^{k'})
		 }.
		 \eqnd
These maps are compatible with thickening, and is natural with respect to families of codimension zero embeddings.
So the equivalence~\eqref{eqn. emb to Fr} is indeed natural.
\end{proof}

\begin{prop}\label{prop. maps out of point preserved}
For every object $X \in \mfld^{\dd}$, the map $\hom_{\mfld^{\dd}}(pt,X) \to \hom_{\Top_{/BO}}(pt,X)$ given by~\eqref{eqn. main functor} is a homotopy equivalence.
\end{prop}

\begin{proof}
Consider the functor~\eqref{eqn. Fr}. If $X$ is modeled by a manifold of dimension $d$, a smooth embedding $j: [0,1]^d \to X$ is sent to the $GL$-equivariant map $\Fr([0,1]^d) \to \Fr(X)$ induced by (the derivative of) $j$. Restriction of $j$ to the origin is compatible with the equivalence~\eqref{eqn. emb to Fr}, in the sense that the following diagram commutes up to homotopy:
		 	\eqnn
			\xymatrix{
			\hom_{\MMfld^{\dd}}(pt,X) \ar[r]^-{\eqref{eqn. emb to Fr}} \ar[d]_{\eqref{eqn. Fr}} & \sing \Fr(X) \\
			\hom_{\TTop^{GL}}(\Fr([0,1]^d), \Fr(X)) \ar[ur]_-{ev_{0,\id_{\RR^d}}}.
			}
			\eqnd
		Here, $ev_{0,\id_{\RR^d}}$ is the evaluation of an $GL$-equivariant map at the canonical frame of the origin of $[0,1]^d$. This evaluation map is a homotopy equivalence, as it is a map of fiber bundles with homotopy equivalent base spaces (equivalent to $X$) and homotopy equivalent fibers (equivalent to $GL$). 
		
		Now we finish by recalling the fact that mapping spaces in the $\infty$-category $\Top_{/BGL} \simeq N(\TTop^{GL})$ may be computed as equivariant mapping spaces between free $GL$-spaces, and that~\eqref{eqn. main functor} is defined by pre-composing with~\eqref{eqn. Fr}. The equivalences~\eqref{eqn. GL spaces to O spaces} and~\eqref{eqn. functor from O-spaces to over BO} finish the job.
		\end{proof}
				
\begin{lemma}\label{lemma. fully faithful}
The functor~\eqref{eqn. main functor} is fully faithful.
\end{lemma}

\begin{proof}[Proof of Lemma~\ref{lemma. fully faithful}.]
Fix smooth manifolds $W$ and $X$. Up to thickening, $W$ is written as a finite colimit of a diagram involving only disks---i.e., only the point, up to thickening. (This is the content of Lemma~\ref{lemma. point generates}.) So write $W \simeq \colim_{\cD} pt$ where $\cD$ is some finite diagram. We have the following homotopy commuting diagram of Kan complexes:
		\eqnn
		\xymatrix{
		\hom_{\mfld^{\dd}}(W, X) \ar[d]
			&\hom_{\mfld^{\dd}}(\colim_{\cD} pt , X) \ar[l]_-{\simeq}  \ar[r]^-{\simeq} \ar[d]
			&\holim_{\cD} \hom_{\mfld^{\dd}}(pt , X) \ar[d] \\
		\hom_{\Top_{/BO}}(W, X) 
			&\hom_{\Top_{/BO}}(\colim_{\cD} pt,X) \ar[l]_-{\simeq}  \ar[r]^-{\simeq} 
			&\holim_{\cD} \hom_{\Top_{/BO}}(pt , X )
		}
		\eqnd
Note that all vertical arrows are obtained from the functor~\eqref{eqn. main functor}. 

The upper-left horizontal arrow is a homotopy equivalence by the hypotheses that $W \simeq \colim_\cD pt$; the two horizontal arrows on the right are equivalences by the definition of colimit. We use Proposition~\ref{prop. main functor preserves finite colims} to conclude the lower-left horizontal arrow is an equivalence. By Proposition~\ref{prop. maps out of point preserved}, the rightmost vertical arrow is a homotopy limit of equivalences---hence an equivalence. It follows that the leftmost vertical arrow is an equivalence, establishing the claim.
\end{proof}

\subsection{The equivalence (proof of Theorem~\ref{theorem. main})}
\label{section. proof of main theorem}
\begin{proof}[Proof of Theorem~\ref{theorem. main}.]
The functor~\eqref{eqn. main functor} 
is fully faithful by Lemma~\ref{lemma. fully faithful}. It remains it is essentially surjective.

By definition, the $\infty$-category of finite spaces
$\Top^{\finite} \subset \Top$
is the full subactegory generated by a point under finite colimits. (For those readers who prefer other models: $\Top^{\finite}$ is equivalent to the $\infty$-category of spaces homotopy equivalent to finite CW complexes. Noting cell attachments are homotopy pushouts along maps $S^{n-1} \to D^n \simeq pt$ and by applying induction on dimensions of cells, we conclude that all finite CW complexes are generated under finite colimits by a point.)

Thus, given an object $X \to BO$ of $\Top^{\finite}_{/BO}$, we have a constant finite diagram $f: D \to \Top^{\finite}$ (with value $pt$) having colimit $X$. Choose a functor $f^{\triangleright} : D^{\triangleright} \to \Top^{\finite}$ exhibiting $X$ as the colimit. Composition with the map $X \to BO$ lifts $f^{\triangleright}$ to a diagram $(\tilde f)^{\triangleright} : D^{\triangleright} \to \Top^{\finite}_{/BO}$. By Remark~\ref{remark. slice categories create colimits}, $(\tilde f)^{\triangleright}$ is a colimit diagram, hence we see that $X \to BO$ is in the subcategory of  $\Top^{\finite}_{/BO}$ generated by the object $pt \to BO$. 

The functor~\eqref{eqn. main functor} contains the point $pt \to BO$ in its image, is fully faithful, and preserves all finite colimits by Proposition~\ref{prop. main functor preserves finite colims}. Thus it is
essentially surjective.
\end{proof}
	 
\bibliographystyle{amsplain} 
\bibliography{thickened-manifolds-biblio}

\end{document}